\documentstyle[11pt,mssymb]{article}

\font \twlmath = msym10 at 11pt
\newfam \mathfam \textfont \mathfam = \twlmath

\newcommand{\arrow}{\longrightarrow}

\newcommand{\Z}{{\Bbb Z}}
\newcommand{\C}{{\Bbb C}}
\newcommand{\R}{{\Bbb R}}
\newcommand{\h}{{\Bbb H}}

\newcommand{\g}{{\goth g}}
\newcommand{\galochka}{\check{\;}}
\newcommand{\6}{\partial}
\newcommand{\1}{\sqrt{-1}\:}

\def\blacksquare{\hbox{\vrule width 4pt height 4pt depth 0pt}}

\begin{document}

\centerline{\bf Hyperholomorphic bundles}
\centerline{\bf over a hyperk\"ahler manifold.}

\hfill

\centerline{Mikhail Verbitsky.
\footnote{ The author is partially supported by NSF grant}}
\centerline{verbit@math.harvard.edu}

\hfill

{\bf 0. Intruduction. }

The main object of this paper is the notion of a hyperholomorphic bundle
(Definition 2.4) over a hyperk\"ahler manifold $M$ (Definition 1.1). The
hyperholomorphic bundle is a direct sum of holomorphic stable
holomorphic bundles. The first Chern class of a hyperholomorphic
bundle is of zero degree.

Roughly speaking, the hyperholomorphic bundle is a bundle which is
holomorphic with respect to all complex structures induced by the
hyperk\"ahler structure on $M$.  It was proven (Proposition
4.1 of [V]) that if M is a complex hyperk\"ahler surface ($K3$ or abelian
surface) then any stable bundle whith first Chern class of zero degree is
hyperholomorphic.

Interesting properties of hyperholomorphic bundles include the analogue of
(p,q)-decomposition, $\6\bar\6$-lemma and an analogue of the
strong Lefshetz theorem on the
holomorphic cohomology $H^*(B)$ of a hyperholomorphic $B$
with a parallel real structure (proven in the Section 4).

For a hyperk\"ahler manifold, one can define
the action of quaternions on its cohomology groups.
The characteristic classes of a hyperholomorphic bundle
are invariant under this action.

Conversely, if $B$ is a stable bundle with the first
two Chern classes invariant under the quaternion
action, the bundle $B$ is hyperholomorphic (Theorem 2.5).

We are describing a coarse moduli space of deformations of a
hyperholomorphic bundle locally (Theorem 6.2). In particular, we show that
there are no obstructions for a deformation besides Yoneda
pairing (Definition 6.2).

This description is used to construct a hyperk\"ahler structure
on the space of deformations of a given hyperholomorphic bundle
(Theorem 6.3), thus generalizing results of [M] and [Ko].

As an application, one can prove that the stable moduli
space for the stable bundles with certain Chern classes
do not depend on the choice of a base manifold in its
deformation class (Proposition 10.3).

Most of results of this paper can be generalized on projectively
hyperholomorphic bundles, which stand in the same relation to
hyperholomorphic ones as projectively flat bundles stand to flat
bundles (Section 11). Over a hyperk\"ahler surface (K3 or
abelian sufrace), any stable Yang-Mills bundle is hyperholomorphic
(Proposition 11.2).

{\bf \centerline {Contents.}}

{\bf 1. Hyperk\"ahler manifolds.

2. Hyperholomorphic bundles.

3. Hermitian bundles: preliminary computations.

4. Hodge analysis and $(p,q)$-decomposition on the
holomorphic cohomology $H^*(B)$ for a hyperholomorphic
$B$.

5. Proof of Theorem 2.5, which states that the stable bundle is
hyperholomorphic if its first two Chern classes are invariant
with respect to the isotropy group of the base hyperk\"ahler manifold.

6. Deformation spaces for holomorphic and hyperholomorphic bundles.

7. Local deformations of a hyperholomorphic connection. The
proof of Theorem 6.2, which states that the only obstruction
to the deformation of a hyperholomorphic bundle is Yoneda product.

8. Comparing Laplacians.

9. The space of deformations of a hyperholomorphic bundle
is hyperk\"ahler (Theorem 6.3).

10. Some applications.

11. Projectively hyperholomorphic bundles.}

\hfill

\hfill

{\bf 1. Hyperk\"ahler manifolds.}

\hfill

{\bf Definition 1.1} ([B], [Bes]) A {\bf hyperk\"ahler manifold} is a
Riemannian manifold $M$ endowed with three complex structures $I$, $J$
and $K$, such that the following holds.

\hspace{5mm}   (i)  $M$ is K\"ahler with respect to these structures and

\hspace{5mm}   (ii) $I$, $J$ and $K$, considered as  endomorphisms
of a real tangent bundle, satisfy the relation
$I\circ J=-J\circ I = K$.

\hfill

This means that the hyperk\"ahler manifold has the natural action of
quaternions ${\Bbb H}$ in its real tangent bundle.
Therefore its complex dimension is even.


Let $\mbox{ad}I$, $\mbox{ad}J$ and $\mbox{ad}K$ be the operators on the
bundles of differential forms over a hyperk\"ahler manifold
$M$ which are defined as follows. Define $\mbox{ad}I$.
Let this operator act as a complex structure operator
$I$ on the bundle of differential 1-forms. We
extend it on $i$-forms for arbitrary $i$ using Leibnitz
formula: $\mbox{ad}I(\alpha\wedge\beta)=\mbox{ad}I(\alpha)\wedge\beta+
\alpha\wedge \mbox{ad}I(\beta)$. Since Leibnitz
formula is true for a commutator in a Lie algebras, one can immediately
obtain the following identities, which follow from the same
identities in ${\Bbb H}$:

\[
   [\mbox{ad}I,\mbox{ad}J]=2\mbox{ad}K;\;
   [\mbox{ad}J,\mbox{ad}K]=2\mbox{ad}I;\;
\]

\[
   [\mbox{ad}K,\mbox{ad}I]=2\mbox{ad}J
\]

Therefore, the operators $\mbox{ad}I,\mbox{ad}J,\mbox{ad}K$
generate a Lie algebra $\goth{su}(2)$ acting on the
bundle of differential forms. We can integrate this
Lie algebra action to the action of a Lie group
$SU(2)$. In particular, operators $I$, $J$
and $K$, which act on differential forms by the formula
$I(\alpha\wedge\beta)=I(\alpha)\wedge I(\beta)$,
belong to this group.

{\bf Proposition 1.1:} There is an action of the Lie group $SU(2)$
and Lie algebra $\goth{su}(2)$ on the bundle of differential
forms over a hyperk\"ahler manifold. This action is
parallel, and therefore it commutes with Laplace operator.
$\blacksquare$

If $M$ is compact, this implies that there is
a canonical $SU(2)$-action on $H^i(M,\R)$ (see [V1]).

\hfill

Let $M$ be a hyperk\"ahler manifold with a Riemannian form $<\cdot,\cdot>$.
Let the form $\omega_I := <I(\cdot),\cdot>$ be the usual K\"ahler
form  which is closed and parallel
(with respect to the connection). Analogously defined
forms $\omega_J$ and $\omega_K$ are
also closed and parallel.

The simple linear algebraic
consideration ([B]) shows that \hfill
$\omega_J+\sqrt{-1}\omega_K$ is of
type $(2,0)$ and, being closed, this form is also holomorphic.
It is called {\bf the canonical holomorphic symplectic form
of a manifold M}. Conversely, if there is a parallel
holomorphic symplectic form on a K\"ahler manifold $M$,
this manifold has a hyperk\"ahler structure ([B]).

If some $compact$ K\"ahler manifold $M$ admits non-degenerate
holomorphic symplectic form $\Omega$, the Calabi-Yau ([Y]) theorem
implies that $M$ is hyperk\"ahler (see [B]).
This follows from the existence of a K\"ahler
metric on $M$ such that $\Omega$ is parallel for the Levi-Civitta
connection associated with this metric.

\hfill

Let $M$ be a hyperk\"ahler manifold with complex structures
$I$, $J$ and $K$. For any real numbers $a$, $b$, $c$
such that $a^2+b^2+c^2=1$ the operator $L:=aI+bJ+cK$ is also
an almost complex structure: $L^2=-1$.
Clearly, $L$ is parallel with respect to a connection.
This implies that $L$ is a complex structure, and
that $M$ is K\"ahler with respect to $L$.

\hfill

{\bf Definition 1.2} If $M$ is a  hyperk\"ahler manifold,
the complex structure $L$ is called {\bf induced
by a hyperk\"ahler structure}, if $L=aI+bJ+cK$ for some
real numbers $a,b,c\:|\:a^2+b^2+c^2=1$.

\hfill

\hfill

If $M$ is a hyperk\"ahler manifold and $L$ is induced complex structure,
we will denote $M$, considered as a complex manifold with respect to
$L$, by $(M,L)$ or, sometimes, by $M_L$.

\hfill

Consider the Lie algebra $\goth{g}_M$ generated by ${ad}L$ for all $L$
induced by a hyperk\"ahler structure on $M$. One can easily see
that $\goth{g}_M=\goth{su}(2)$.
The Lie algebra $\goth{g}_M$ is called {\bf isotropy algebra} of $M$, and
corresponding Lie group $G_M$ is called an {\bf isotropy group}
of $M$. By Proposition 1.1, the action of the group is parallel,
and therefore it commutes with Laplace operator in differential
forms. In particular, this implies that the action of the isotropy
group $G_M$ preserves harmonic forms, and therefore this
group canonically acts on cohomology of $M$.

\hfill

{\bf Proposition 1.2:} Let $\omega$ be a differential form over
a hyperk\"ahler manifold $M$. The form $\omega$ is $G_M$-invariant
if and only if it is of Hodge type $(p,p)$ with respect to all
induced complex structures on $M$.

{\bf Proof:} Assume that $\omega$ is $G_M$-invariant.
This implies that all elements of $\g_M$ act trivially on
$\omega$ and, in particular, that $\mbox{ad}L(\omega)=0$
for any induced complex structure $L$. On the other hand,
$\mbox{ad}L(\omega)=(p-q)\1$ if $\omega$ is of Hodge type $(p,q)$.
Therefore $\omega$ is of Hodge type $(p,p)$ with respect to any
induced complex structure $L$.

Conversely, assume that $\omega$ is of type $(p,p)$ with respect
to all induced $L$. Then $\mbox{ad}L(\omega)=0$ for any induced $L$.
By definition, $\g_M$ is generated by such $\mbox{ad}L(\omega)=0$,
and therefore $\g_M$ and $G_M$ act trivially on $\omega$. $\blacksquare$

\hfill

\hfill

{\bf 2. Hyperholomorphic bundles.}

\hfill

 Let $B$ be a Hermitian vector bundle over the complex
manifold $M$. Let $\theta$ be a Hermitian connection
on $B$ and $\Theta\in\Lambda^2\otimes End(B)$ be its curvature.
If $B$ is holomorphic, this connection
is called {\bf compatible with a holomorphic structure} if
$\nabla_X(\zeta)=0$ for any holomorphic section $\zeta$ and
any antiholomorphic tangent vector $X$. If there exist
a holomorphic structure compatible with the given
Hermitian connection then this connection is called
{\bf integrable}.

\hfill

One can define a {\bf Hodge decomposition} in the space of differential
forms with coefficients in any complex bundle, in particular, $End(B)$
(see [GH]).

{\bf Theorem 2.1} (Newlander-Nirenberg)
The connection $\theta$ in $B$ is integrable
if and only if $\Theta\in\Lambda^{1,1}\otimes End(B)$, where
$\Lambda^{1,1}\otimes End(B)$ denotes (1,1)-forms with
respect to the Hodge decomposition.
The holomorphic structure compatible with a
given connection $\theta$ is unique.

{\bf Proof:} This is Proposition 4.17 of [Ko], Chapter I.
$\blacksquare$

\hfill

\hfill

{\bf Definition 2.1} Let $B$ be a Hermitian bundle with
a Hermitian connection $\theta$ over a hyperk\"ahler manifold
$M$. The connection $\theta$ is called {\bf hyperholomorphic} if it is
integrable with respect to any complex structure induced
by a hyperk\"ahler structure.

As follows from the Theorem 2.1, $\theta$ is hyperholomorphic
if and only if its curvature $\Theta$ is of type (1,1) with
respect to any of complex structures induced by a hyperk\"ahler
structure.

As follows from Proposition 1.2, $\theta$ is hyperholomorphich
if and only if $\Theta$ is a $G_M$-invariant differential form.

\hfill

Let $M$ be a K\"ahler manifold with a K\"ahler form $\omega$.
For differential forms with coefficients in any vector bundle
there is Hodge operator $L: \eta\arrow\omega\wedge\eta$.
There is also fiberwise-adjoint Hodge operator $\Lambda$
(see [GH]).

\hfill

{\bf Definition 2.2} (see [UY]) Let $B$ be a holomorphic bundle over $M$
with a holomorphic Hermitian connection $\theta$ and a
curvature $\Theta\in\Lambda^{1,1}\otimes End(B)$.
The Hermitian metric on $B$ and the connection $\theta$
defined by this metric is called Yang-Mills if

\[\Lambda(\Theta)=constant\cdot Id,\]

\hspace{-6mm}where $\Lambda$ is a Hodge operator and $Id$ is
the identity endomorphism which is a section of $End(B)$.

\hfill

{\bf Definition 2.3}  Let $F$ be a coherent sheaf over
$n$-dimensional compact K\"ahler manifold $M$. We define
$deg(B)$ as

\[ deg(B)=\int_M\frac{ c_1(F)\wedge\omega^{n-1}}{vol(M)}\]

and $slope(B)$ as
\[ slope(B)=\frac{1}{rank(F)}\cdot deg(B) \]

The number $slope(F)$ is rational, and it depends only on a
cohomology class of $c_1(F)$. The simple bundle $B$ is called {\bf stable}
if for every coherent subsheaf $B'\subset B$ with $rank(B') < rank(B)$
\[slope(B')<slope(B)\]

\hfill

Later on, we will usually consider the bundles $B$ with $deg(B)=0$.

{\bf Proposition 2.1:} Let $M$ be a compact K\"ahler manifold,
and $B$ be a Hermitian holomorphic bundle over $M$
with $deg(B)=0$. The bundle $B$ is Yang-Mills if and only if
$\Lambda(\Theta)=0$.

{\bf Proof:} One can see this
using the following argument. The number

\[
   deg(B)=\int_M \frac {c_1(B)\wedge\omega^{n-1}}{vol(M)}
\]
is equal by the Gauss-Bonnet theorem to the value of an integral
\[
   \int_M \frac{\Lambda(Tr(\Theta))}{Vol(M)}.
\]
For a Yang-Mills bundle, the equation $deg(B)=0$ implies
that $\int_M \Lambda(\Theta) = 0$. $\blacksquare$

\hfill

A holomorphic bundle is called  {\bf undecomposable}
if it cannot be decomposed onto a direct sum
of two or more holomorphic bundles.

There is an important theorem ([UY])

{\bf Theorem 2.2} (Uhlenbeck-Yau): Let B be some undecomposable
holomorphic bundle over a compact K\"ahler manifold. Then $B$ admits
a Yang-Mills Hermitian metric if and only if it is stable, and
this metric is unique (up to a constant).

\hfill

{\bf Theorem 2.3} A hyperholomorphic connection in a bundle B
is Yang-Mills. Moreover, for such connection $\Lambda(\Theta)=0$
where $\Theta$ is a curvature in $B$.

{\bf Proof:} We will use the definition of a hyperholomorphic
connection as one with $G_M$-invariant curvature.
Theorem 2.3 follows from the

{\bf Lemma 2.1 } Let $\Theta$ be a $G_M$-invariant
section of $\Lambda^2(M)\otimes End(B)$. Then
$\Lambda_L(\Theta)=0$ for each induced complex structure
$L$. By $\Lambda_L$ we understand the Hodge operator $\Lambda$
associated with the K\"ahler complex structure $L$.

{\bf Proof of Lemma 2.1} The action of $G_M$, as defined above, is
orthogonal with respect to the standard metric on the space of
differential forms. From the theory of group representations we know
the following fact. Let $G$ be a group which acts orthogonally on the
Hermitian space $V$. Let $V_{inv}$ be the space of $G$-invariants
in $V$ and $V=V_{inv}\oplus V_1$ be a $G$-invariant decomposition. Then
$V_1$ is orthogonal to $V_{inv}$. We will use this to prove
Lemma 1.2.

Consider the bundle $R\subset \Lambda^2(M)\otimes End(B)$ spanned
by the K\"ahler forms  $\omega_L$, for all induced complex structures $L$
on $M$. The element of a
group $G_M$ maps each of these forms into a linear combination
of these forms. This could be proven by a direct computation.
Therefore, $R$ is invariant with respect to the $G_M$-action.
Moreover, no section of $R$ is $G_M$-invariant. This is
proven by another direct computation. Both of these computations
are based on the fact that $R$ is a trivial bundle
over $M$ with a fiber $\g_M$, and $G_M$ acts on $R$ by means of
adjoint representation on $\g_M$.

By the statement in the beginning of the proof,
this bundle is orthogonal to  the bundle of $G_M$-invariants.

Now, by the definition of the Hodge operator $\Lambda_L$,
for $G_M$-invariant $\Theta$,

\[ \Lambda_L(\Theta)=0 \]
because the K\"ahler form $\omega_L$ is
orthogonal to $\Theta$.
Theorem 2.3 is proven. $\;\;\blacksquare$

\hfill

{}From now on the base hyperk\"ahler manifold $M$ is supposed to be
compact.

\hfill

{\bf Theorem 2.4} Let $M$ be a hyperk\"ahler surface, i. e., $K3$ or
abelian surface. A stable holomorphic bundle B with $deg(B)=0$
over $M$ always admits a hyperholomorphic connection,
which is unique. ([V], Theorem 2.1)

{\bf Proof:}

Take Yang-Mills Hermitian connection on $B$. We shall prove that
it is hyperholomorphic. The uniqueness of a hyperholomorphic connection
follows from Theorem 2.3, because it is Yang-Mills, and Yang-Mills
Hermitianconnection in a stable bundle is unique.

The space of complex 2-forms over $M$
has a trivial 3-dimensional subbundle $P$
spanned by K\"ahler form, holomorphic
symplectic form and the form conjugated to the holomporhic
symplectic form. Let $P^\bot$ be its orthogonal completion.
One easily sees that if $M$ is a surface then $P^\bot$
is a bundle of $G_M$-invariant 2-forms. On the other hand,
the curvature $\Theta$ of a Yang-Mills bundle of degree zero is
orthogonal to $P$ for any hyperk\"ahler manifold
by the following reasons. First,
$\Theta$ is of type (1,1) and this is sufficient for
the orthogonality to the holomorphic symplectic form and its
conjugate. Second, $\Lambda(\Theta)=0$ is equivalent
to $\Theta$ being orthogonal to the K\"ahler form by definition
of $\Lambda$. $\blacksquare$.

\hfill

One should note that in the proof above, $P$ is identical to
$\Lambda^+(M)$ and $P^\bot$ is $\Lambda^-(M)$, where
$\Lambda^2(M)=\Lambda^+(M)\oplus \Lambda^-(M)$  is a standard
decomposition of 2-forms over a Riemannian 4-fold on
autodual and anti-autodual forms.

\hfill

For any stable holomorphic bundle there exist unique
HermitianYang-Mills connection which, for some bundles,
turns out to be hyperholomorphic. This observation gives us
independent from connection notion
of a stable hyperholomorphic bundle.

\hfill

Let M be a hyperk\"ahler manifold, and let $I$ be an
induced complex structure over $M$.

{\bf Definition 2.4} The stable holomorphic bundle $B$
over $(M,I)$ is called {\bf simple hyperholomorphic} if the unique
Yang-Mills connection
on $B$ is hyperholomorphic. Generally, $B$  is called {\bf
hyperholomorphic} if it could be decomposed
onto a direct sum of simple hyperholomorphic stable bundles.

\hfill

Further on, we will consider mostly simple hyperholomorphic bundles
and omit the word ``simple'' when it is obvious from context.

\hfill

For a hyperholomorphic bundle, its $p$-th Chern class is of type
$(p,p)$ with respect to any of induced complex structures.
By Proposition 1.2, this implies that Chern classes of
a hyperholomorphic bundle are $G_M$-invariant. Conversely,
there is a theorem of extreme importance, which
is proven in the Section 5:

{\bf Theorem 2.5:} Let $B$ be a stable bundle over
$(M,I)$, where $M$ is a hyperk\"ahler manifold and $I$
is an induced complex structure over $M$. If $c_1(B)$
and $c_2(B)$ are $G_M$-invariant, the stable bundle $B$
is hyperholomorphic.

\hfill

One knows that $c_i(B)$ are $G_M$-invariant for any hyperholomorphic
bundle (see the paragraph above the theorem). One also knows
that the form is $G_M$-invariant if and only if it is
of Hodge type $(p,p)$ with respect to all induced complex
structures (Proposition 1.2).

\hfill

\hfill

{\bf 3. Hermitian bundles: preliminary observations.}


Throughout this section, $M$ is a K\"ahler manifold
and $E$ is a Hermitian holomorphic bundle over $M$ with a Hermitian
holomorphic connection $\nabla$ and curvature $\Theta$. Consider the
Hodge decomposition for $\nabla$ ([GH])

(3.1)\hspace{5mm} $\nabla = \6 + \bar\6$

(traditionally, the holomorphic component of $\nabla$ is denoted by
$\nabla'$ while antiholomorphic component is denoted by $\bar\6$).  We will
denote the fiberwise-adjoint operators by $\6^*$ and $\bar\6^*$.

This notation is different from a traditional one.

\hfill

Our major tool are famous Kodaira's identities (see [GH]):

\[
    [ \Lambda, \6] = \1\bar\6^* {\rm \;\; and \;\;}
    [ \Lambda, \bar\6] = -\sqrt{-1}\6^*
\]

where $\Lambda$ is the Hodge operator.

\hfill

Consider four Laplace operators:

\[ \Delta_\6 = \6\6^*+\6^*\6; \; \;
\Delta_{\bar\6} = \bar\6\bar\6^* + \bar\6^*\bar\6; \; \;
\Delta_d = \nabla\nabla^* + \nabla^*\nabla;
\]
\[
  \Delta_{d^c}= I\circ \Delta_d \circ I^{-1}
\]
where $I$ is a complex structure operator acting
on differential forms by multiplicativity.

If $B$ is trivial, $\Delta_\6=\Delta_{\bar\6}=1/2\;\Delta_d=
1/2\:\Delta_{d^c}$ (see [GH]). In case of nontrivial
Hermitian holomorphic bundle $B$, Kodaira identities
imply the following.

\hfill

{\bf Proposition 3.1}

\hspace{5mm}a)  $\Delta_\6+\Delta_{\bar\6}=\Delta_d$

\hspace{5mm}b)   $\Delta_\6-\Delta_{\bar\6}=\1[\Lambda;\nabla^2]$

\hspace{5mm}c)  $\Delta_d=\Delta_{d^c}$,

\hspace{-3mm}where $\Lambda$ is a Hodge operator.

{\bf Proof}

a) By (3.1), $\Delta_d=\Delta_\6+\Delta_{\bar\6} + (\6\bar\6^*+\bar\6^*\6)+
(\6^*\bar\6+\bar\6\6^*)$ We need only to prove that the last two terms are
zero.  By the Kodaira identity,

\[\Lambda\6-\6\Lambda= \sqrt{-1}\bar\6^*\]

and therefore

\[\frac{\6\bar\6^*+\bar\6^*\6}{\sqrt{-1}}=
\6\Lambda\6- \6^2\Lambda + \Lambda\6^2 -\6\Lambda\6 = 0\]

Vanishing of the last term is analogous.

\hfill

b) By Kodaira identities, $\;\1\Delta_{\bar\6} = \bar\6(\Lambda\6-\6\Lambda)
+(\Lambda\6-\6\Lambda)\bar\6$,

while $\;-\1\Delta_\6 = \6(\Lambda\bar\6-\bar\6\Lambda)
+(\Lambda\bar\6-\bar\6\Lambda)\6$.

The summation of these equations yields

\[ \1(\Delta_{\bar\6}-\Delta_{\6})= - (\bar\6\6 + \6\bar\6)\Lambda +
\Lambda(\6\bar\6 + \bar\6\6) = [\Lambda;\nabla^2].\]

\hfill

c) Let us take $d^c=I\circ \nabla\circ I^{-1}$.
For symmetry, denote $\nabla$ by $d$. Kodaira identities state
that $[\Lambda,  d ]= {d^c}^*$ and $[\Lambda, d^c]= - d ^*$.
Analogously to the case b),
\[
   \Delta_{d^c} = d^c(\Lambda d - d \Lambda)
   +(\Lambda d - d \Lambda)d^c,
\]
while
\[
   -\Delta_d= d (\Lambda d^c-d^c\Lambda)
   +(\Lambda d^c-d^c \Lambda) d.
\]
Summing up, we obtain
$\Delta_{d^c} - \Delta_d= [\Lambda, d^c d + d d^c]$.
Clearly,

\[ d^c d+d d^c=\frac{\6-\bar\6}{\1}(\6+\bar\6)
+(\6+\bar\6)\frac{\6-\bar\6}{\1}= -\1(\6\bar\6-\bar\6\6+
\bar\6\6-\6\bar\6)=0
			\;\;\blacksquare\]

\hfill

One knows that $\nabla^2$ is a linear operator on forms which
maps the form $\eta$ into $\Theta\wedge\eta$.
Therefore the Proposition 3.1 implies that the difference
between the Laplacians $\nabla_\6$, $\nabla_{\bar\6}$ and
$\nabla_d$ is a differential operator of degree zero,
i. e. the linear operator. We will use analogous computations
in the section 4.

\hfill

\hfill

The rest of this section is dedicated to the real structures
in Hermitian bundles.

\hfill

The {\bf real structure} on the complex vector bundle is
anticomplex involution. In other words, the real
structure on a complex bundle $B$ is an operator $T$
such that $T^2=Id$ and $T(\lambda k)=\bar\lambda x$
for $\lambda\in\C$.

Suppose that $B$ is a bundle with connection. We
will call a real structure $T$ {\bf parallel} if $T$ is
parallel as a section of $End_\R(B,B)$ with respect
to a connection in $End_\R(B,B)$ induced by a
connection in $B$.

For a real structure $T$ on a bundle $B$, let
$B_{\Bbb R} := \{b\in B\;|\; T(b)=b\}$ be the space of $T$-invariant
vectors. Clearly, $B_\R$ is endowed
with a unique connection such that $B=B_{\Bbb R}\otimes \C$
as a bundle with connection.

For a real structure $T$, the operator $-T$ also defines a
real structure.

If $E$ is a Hermitian bundle, the bundle $End(E)$ is
endowed with a canonical real structure $T$. This real structure
maps an endomorphism $\alpha$ of $E$ to an endomorphism
$-\alpha^\bot$. By $\alpha^\bot$ we denote the endomorphism
which is adjoined to $\alpha$ with respect to
the Hermitian product on $E$.

Consider the connection in $End(E)$ which is induced by a
connection in $E$. The constructed above real structure operator is
parallel with respect to this connection.

One should note that the connection in $End(E)$ is Yang-Mills
if it is induced by a Yang-Mills connection in $E$.

Generally speaking, one can easily tell whether a
Yang-Mills bundle $B$ admits a parallel real structure:

{\bf Proposition 3.2} Let $B$ be a Yang-Mills bundle
over a K\"ahler manifold $M$. The bundle $B$ is endowed with
a parallel real structure if and only if $B$, taken as a
holomorphic bundle, is isomorphic to its dual bundle,
and this isomorphism is defined by a symmetric
non-degenerate bilinear form.
If $B$ is nondecomposable, this real structure
is unique up to a multiplication by $-1$.

{\bf Proof:} The bundle $B$ is always canonically
isomorphic to $B^*$ as a real bundle with connection. The
isomorphism is given by the real part of the
Hermitian form on $B$. Let $f:B\arrow B^*$
be the isomorphism operator. Clearly, $I\circ f=-f\circ I^*$,
where $I$ denotes the complex structure operator
(multiplication by $\1$)
on $B$ and $I^*$ denotes the complex structure operator
 on $B^*$. Now, if $B$ has a parallel real structure $T$,
the operator $T\circ f$ commute with the complex structure operator
because
\[
   I\circ T\circ f=- T\circ I\circ f= T\circ f\circ I^*.
\]
Since this $T\circ f$ is parallel, it defines an isomorphism
of Hermitian bundles $B$ and $B^*$.
The converse claim follows from the uniqueness
of Yang-Mills metric. $\;\;\blacksquare$

\hfill

If $B$ is endowed with a parallel real structure,
the bundle $\Lambda^*(M,B)$
of differential forms with coefficients in $B$
is also endowed with a parallel real structure.
It could be most easily seen on the following way.
Let $\Lambda^i(M,\R)$ be the bundle
of real differential $i$-forms over $M$. If $B$ has
a real structure, $B=B_\R\otimes \C$. Therefore,
\[
   \Lambda^i(M,\R)\otimes_\R B_\R\otimes_\R\C
   = \Lambda^i(M,B).
\]

Consider the operator $T$ which acts as a complex conjugation
in $\C$ and trivially on the first two components of a tensor
product $\Lambda^i(M,\R)\otimes B_\R\otimes_\R \C$.

Obviuosly, the operator $T$ defines a parallel real structure on
$\Lambda^i(M,B)$.

One can easily prove that $T$ maps $\Lambda^{pq}(M,B)$ into
$\Lambda^{qp}(M,B)$.

\hfill

\hfill

{\bf 4. Hodge analysis and (p,q)-decomposition on the holomorphic
cohomology $H^*(B)$ for a hyperholomorphic $B$.}

\hfill

In this section we will prove several theorems about differential\\
(p,0)-forms with coefficients in a hyperholomorphic bundle $B$.  These
theorems are direct analogues of
Kodaira identities, strong Lefshetz theorem,
(p,q)-decomposition and $dd^c$-lemma for forms with trivial coefficients
over a K\"ahler manifold (see [GH]).

The difference between the traditional and our cases is following.
Traditionally, these theorems are being proven for the
space of (topological) cohomology $H^*(M,\C)$
of a compact K\"ahler manifold. In our case, their analogues
exist for the space $H^*(M,B)$ of holomorphic cohomology
of $B$.

\hfill

Let $B$ be a holomorphic Hermitian bundle with a connection and
a parallel real structure over a hyperk\"ahler manifold
$(M,I)$.  Let $I$, $J$ and $K$ be complex structure operators in a real
cotangent bundle over $M$, as usually, $I\circ J=-J\circ I=K$.
These operators could be extended on the bundles of differential
forms over $M$ by multiplicativity, i. e., by the formulas
$J(\alpha\wedge\beta)=J(\alpha)\wedge J(\beta)$ etc..  Since
$I\circ J=-J\circ I$, the operator $J$ maps $(1,0)$-forms on $(0,1)$-forms,
and therefore it maps $(p,0)$-forms on $(0,p)$-forms.

Consider the operator
\[\bar J:\; \Lambda^{(p,0)}(M,B)\arrow \Lambda^{(p,0)}(M,B)\]
which is a composition of $J$ and a real structure operator.
Let
\[ \6^j=\bar J\circ\6\circ\bar J^{-1},\]
where $\6:\;\Lambda^{(p,0)}(M,B)\arrow\Lambda^{(p+1,0)}(M,B)$ is a
(1,0)-component of a connection.

{\bf Proposition 4.1} The following equations hold:

\[ (\6^j)^2=0\]
\[\6^j\6+\6\6^j=0\]

{\bf Proof:} The first formula is an obvious consequence
of the relation $\6^2=0$. Let
$\Lambda\Lambda^p(B)=\Lambda(M,B)\otimes_\R\C$. Since $\bar J^2=-1$, the
bundle $\Lambda\Lambda^1(B)$ can be decomposed onto the direct sum

\[\Lambda\Lambda^1(B)=\Lambda\Lambda^{1,0}(B)\oplus\Lambda\Lambda^{0,1}(B),\]

where $\Lambda\Lambda^{1,0}(B)$ corresponds to the eigenspace of an
eigenvalue $\1$ for $\bar J$ and $\Lambda\Lambda^{0,1}(B)$
to eigenvalue $(-\1)$. Analogously, let us by multiplicativity
decompose

\[\Lambda\Lambda^n(B)=\oplus_{p+q=n}\Lambda\Lambda^{p,q}(B)\]
like it is being done in a usual Hodge theory.

Let $\delta=\frac{\6+\1\6^j}{2}$, and $\bar\delta=\frac{\6-\1\6^j}{2}$ be
the (1,0) and (0,1) components of $\6$ with respect to
$(p,q)$-decomposition on $\Lambda\Lambda^*(B)$.
By definition, $\delta^2$ is $(2,0)$-component of $\6^2$
with respect to the $\Lambda\Lambda$-decomposition.
Analogously, $\bar\delta^2$ is $(0,2)$-component
of $\6^2$. Since $\6^2=0$, we see that
$\delta^2=\bar\delta^2=0$. On the other hand,

\[\delta^2= \6^2+(\6^j)^2+\1\6\6^j+\1\6^j\6 \]
and therefore

\[\6^j\6+\6\6^j=0.\;\;\blacksquare\]

Let $L_I$ be the usual Hodge operator $L$ in the space of differential
forms $\Lambda^*(M,B)$ (exterrior multiplication by the K\"ahler form
$\omega_I$), $L_J$ and $L_K$ be analogous operators associated with
complex structures $J$ and $K$ on $M$.

{\bf Proposition 4.2} The following equation holds:

\[ [L_J,\6^*]=\6^j.\]

{\bf Proof:} Let $(\;)^J$ denote $J\circ(\;)\circ J^{-1}$ where $(\;)$ is
some operator on differential forms, $(\;)^I$ and $(\;)^K$ be
analogously defined operators.
Let  $\nabla$ denote a connection in $B$,
and $\nabla^*$ denote the adjoint operator. By
definition, $\6=\frac{\nabla+\1 \nabla^I}{2}$, and therefore
$\6^j=\frac{\nabla^J-\1 \nabla^K}{2}$. On the other
hand,

\[
   [L_J, \6^*]= \left[ L_J, \frac{\nabla^*+\1 (\nabla^*)^I}{2}\right]
   = \frac{\nabla^J + \1 [L_J^{I^{-1}}, \nabla^*]^I}{2}
\]
because by Kodaira identities
$[L_R,\nabla^*]=\nabla^R$ for $R= I, J$ or $K$.
We can apply Kodaira's identity because $B$ is hyperholomorphic
and therefore $\nabla$ is a holomorphic Hermitian connection
with respect to $I$, $J$ and $K$. Now, the simple calculation
shows that $L_J^{I^{-1}}=-L_J$, and therefore

\[
   [L_J,\6^*]= \frac {\nabla^J+\1 \nabla^{J\circ I}}{2}
   =\frac {\nabla^J-\1 \nabla^K}{2}. \;\;\;\blacksquare
\]

{\bf Corollary 4.1} The following is also true:

\[ [L_J,\delta^*]=\1\bar\delta\]
\[ [L_J,\bar\delta^*]=-\1\delta.\;\;\;\blacksquare\]

\hfill

{\bf Proposition 4.3}

a). The operators $\6^*$ and $\6^j$ anticommute, i. e.,
$\6^*\6^j+\6^j\6^*=0$. Analogously, the pairs $(\6^j)^*$ and $\6$,
$\delta^*$ and $\bar\delta$, $\bar\delta^*$ and $\delta$ anticommute.

b) Operators $\6$, $\6^j$, $\delta$ and $\bar\delta$
pairwise anticommute.

{\bf Proof} a) For instance, let us prove anticommutation of the operators
in the first pair. By Proposition 4.2,

\[ \6^*\6^j+\6^j\6^*= \6^* L_J \6^* - \6^*\6^*L_J +
L_J\6^*\6^* - \6^*L_J\6^*=0 \]

b). The operator $\delta\bar\delta+\bar\delta\delta = 1/2\cdot
(\6\6^j+\6^j\6)$ is just (1,1) component of an operator $\6^2$ for a
$(p,q)$-decomposition on $\Lambda\Lambda(B)$ (see proof of Proposition 4.1).
Since the connection in $B$ is holomorphic, $\6^2=0$. Therefore the $(1,1)$
component of $\6^2$, namely $\delta\bar\delta+\bar\delta\delta$, is also zero.
The anticommutation of other pairs is analogous. $\;\;\;\blacksquare$

\hfill

Let $\Delta_{\6^j}$ be the Laplacian associated with $\6^j$.  Let
$\Delta_\delta$ and $\Delta_{\bar\delta}$ be Laplacians on the complex
$\Lambda \Lambda^*(B)$, associated with
$\delta$ and $\bar\delta$ respectively. The operators $\Delta_\6$ and
$\Delta_{\6^{j}}$ could also be considered as operators on
$\Lambda\Lambda^*(B)=\Lambda^{p,0}(M,B)\otimes_\R\C$.

{\bf Theorem 4.1} The following Laplace operators are equal
(or proportional):

\[\Delta_{\6^j}=\Delta_\6=2\Delta_\delta=2\Delta_{\bar\delta}\]

{\bf Proof}: By definition of $\delta$ and $\bar\delta$,

\[ \Delta_\6=\Delta_\delta+\Delta_{\bar\delta} +
      (\delta\bar\delta^*+\bar\delta^*\delta)+
(\delta^*\bar\delta+\bar\delta\delta^*) = \Delta_\delta+\Delta_{\bar\delta}
\] since, by the Proposition 4.3a,
$(\delta\bar\delta^*+\bar\delta^*\delta)=
(\delta^*\bar\delta+\bar\delta\delta^*) = 0$
\hfill

By Corollary 4.1, $\;\1\Delta_{\bar\delta} =
\bar\delta^*(L_J\delta^*-\delta^*L_J)
+(L_J\delta^*-\delta^*L_J)\bar\delta^*$, \\ while $\;-\1\Delta_\delta^* =
\delta^*(\Lambda_J\bar\delta^*-\bar\delta^*\Lambda)
+(L_J\bar\delta^*-\bar\delta^*L_J)\delta^*$.

The summation of these equations yields

\[ \1(\Delta_{\bar\delta}-\Delta_\delta)=
- (\bar\delta^*\delta^* + \delta^*\bar\delta^*)L_J+
L_J(\delta^*\bar\delta^* + \bar\delta^*\delta^*) = 0,\] since
$(\delta^*\bar\delta^*+\bar\delta^*\delta^*)^*$ is equal zero by the
Proposition 4.3(b). Therefore, $\Delta_{\bar\delta}=\Delta_\delta$
and $\Delta_\6=2\Delta_\delta$. It remains only to prove that
$\Delta_\6=\Delta_{\6^j}$. This follows from equations

\[
   \Delta_\6=\6\6^*+\6^*\6=\6[L_J\6^j]+[L_J\6^j]\6=
   \6L_J\6^j-\6^jL_J\6
\]
and
\[
   \Delta_{\6^j}=\6^j{\6^j}^*+{\6^j}^*\6^j=-\6^j[L_J\6]-[L_J\6]\6^j=
   \6L_J\6^j-\6^jL_J\6. \;\;\blacksquare
\]

\hfill

{\bf Corollary 4.2} (The $(p,q)$-decomposition on cohomologies.) The
cohomology space $H^n(B)\otimes_\R\C$ admits $(p,q)$-decomposition:

\[H^n(B)\otimes_\R\C=: \oplus_{p+q=n}H^{p,q}_{hyp}(B).\]

{\bf Proof} Laplacian $\Delta_\delta$ preserves $(p,q)$-decomposition
defined on $\Lambda\Lambda^n(B)$, but its kernel is $H^n(B)\otimes_R\C$ by
the Theorem 4.1. $\;\;\;\blacksquare$

\hfill

{\bf Corollary 4.3} Let $B$ be a hyperholomorphic
bundle with a parallel real structure.
There is a canonical action of the Lie
group $SU(2)$ on the holomorphic cohomology $H^*(B)$ of $B$ as
on $\R$-linear space.

{\bf Proof:} Consider the action of operators $\bar J$ and $I$
on $\Lambda^{p,0}(M,B)$. Clearly, $I\circ\bar J= -\bar J\circ I$.
Let $\bar K := I\circ\bar J$. One can easily check that the
set of operators
\[
   \{aI+b\bar J+c\bar K+d\cdot\mbox{\it Id}\;|\; a^2+b^2+c^2+d^2=1 \}
\]
form a group isomorphic to $SU(2)$. By Theorem 4.1 the action
of this group maps $\Delta_\6$-harmonic forms from
$\Lambda^{p,0}(M,B)$ into $\Delta_\6$-harmonic ones.
This defines the action of $SU(2)$ on the space
of $\Delta_\6$-harmonic forms from $\Lambda^{p,0}(M,B)$,
which is, by Hodge theory ([GH]), canonically isomorphic
to the space of holomorphic cohomology $H^p(B)$.
$\;\;\blacksquare$

One notes that the $SU(2)$-action on $H^1(B)$ extends
to a quaternion action.

\hfill

Let $L_c=L_J+\1 L_K$ be an operator of type (2,0) (see section 1) mapping
the $\Lambda^{p,0}(M,B)$ into $\Lambda^{p+2,0}(M,B)$,
and $\Lambda_c = \Lambda_J+\1\Lambda_K$ be
its adjoint operator.  Proposition 4.2 together with the Theorem 4.1 shows
that $\Lambda_c$ preserves harmonicity of forms, and therefore we can
correctly define the action of $L_c$ and $\Lambda_c$ on the space of
holomorphic cohomologies $H^*(B)$. A simple computation shows that

\[H=[L_c,\Lambda_c]\;\;|\;\;_{H^i(B)} = (n-2i),\]
where $(n-2i)$ means multiplication by a scalar $n-2i$, and $n$ is a
dimension of $M$. There is an Hermitian metric on $H^i(M,B)$
as on the space of $\Delta_\6$-harmonic sections of
the Hermitian bundle $\Lambda^{i,0}(M,B)$.

{\bf Theorem 4.2} (analog of a strong Lefshetz theorem):
Thus defined operators $L_c$, $\Lambda_c$ and $H$
generate a Lie algebra ${\goth sl}(2)$ which acts on $H^*(M,B)$
The standard Hodge-Riemann relations between the metric,
Hodge decomposition and $\goth{sl}(2)$-action also hold.

{\bf Proof} The proof is completely
analogous to the proof of Lefshetz theorem
and Hodge-Riemann relations, see [GH].$\;\;\blacksquare$

The generalization of this theorem could be proven for any
hyperholomorphic bundle $B$ regardless of the existence of
a parallel real structure. Consider the obvious map
$E\otimes {\cal O}_M\arrow E$, where $E$ is a sheaf over $M$ and
${\cal O}_M$ is a structure sheaf. This map defines a
K\"unneth map of cohomology

\[ H^i(E)\times H^j({\cal O}_M)\arrow H^{i+j}(E). \]

One knows that for a hyperk\"ahler manifold of complex
dimension $2n$ the ring $H^*({\cal O}_M)$
contains the sub ring of truncated polynomials of
one variable:

\[ H^*({\cal O}_M)=\C[x]/\{x^{n+1}=0\}. \]
There, $x$ is a generator of $H^2({\cal O}_M)$
corresponding to a canonical holomorphic symplectic form over
$M$ (see [B]). Consider the multiplication by $x$ as a map
$L_c:\; H^i(E)\arrow H^{i+2}(E)$. This map is defined for
any sheaf $E$ over a hyperk\"ahler manifold $M$.

{\bf Theorem 4.2.A} Let $E$ be a hyperholomorphic bundle.
For $i\leq n$ the map

\[ L_c:\; H^i(E)\arrow H^{i+2}(E) \]
is injection. The map

\[ L_c^{n-i}:\; H^i(E)\arrow H^{2n-i}(E) \]
is isomorphism.

{\bf Proof:} Let $B= E\oplus E^*$, where $E^*$ is a dual bundle.
The bundle $B$ has a canonical real structure as following.
Consider the Hermitian form as the anticomplex operator
from $E$ to $E^*$. Denote this operator by $\tau$.
The real structure $T$ in $B$ is defined by the formula
\[ T(\alpha,\beta)= (\tau^{-1}(\beta),\tau(\alpha)) \]

Now, applying Theorem 4.2 to the bundle $B$, we obtain
the statement of Theorem 4.2.A.

\hfill

{\bf Theorem 4.3} ($\6\6^j$-lemma) Let $\omega\in\Lambda^{p,0}(M,B)$ be
some $\6$- and $\6^j$-closed form. Suppose that $\omega$ is also either
$\6$- or $\6^j$-exact. Then there is a form $\kappa$ such that
$\6\6^j(\kappa)=\omega$.

{\bf Proof}.

One should keep in mind that by virtue of Proposition 4.3(a),
$\6\6^j(\eta)=-\6^i\6(\eta)$ for any differential form
$\eta\in\Lambda^i(M,B)$. By the same argument, the Laplacian
$\Delta_\6$ commutes with $\6$ and $\6^j$.

Let $G_\6$ be the Green operator associated with the Laplacian
$\Delta_\6$, and $G_{\6^j}$ be one associated with $\Delta_{\6^j}$.  By
Theorem 4.1, $G_{\6}=G_{\6^j}$. Let us denote this operator by $G$.
Obviously, $G$ commutes with $\6$ and $\6^j$.  Since
$\omega$ is $\6$- (or $\6^j$-) exact, $G\Delta\omega=\omega$, and therefore
its orthogonal projection on the space of harmonic forms is zero. This
proves that $\omega$ is both $\6$- and $\6^j$-exact.

On the other hand, $\Delta_\6 \omega = \6\6^*\omega$ because $\omega$ is
closed. Therefore, $\omega = G\Delta_\6\omega = \6\6^* G\omega$. Let $\omega
= \6^j \gamma$ for some $\gamma$. Since $\6^j$ anticommutes with $\6$ and
$\6^*$ (Proposition 4.3),

\[\omega=G\Delta\omega=\6\6^* G\omega =
\6\6^* G \6^j\gamma = - \6\6^j\6^* G\gamma. \;\;\blacksquare \]

\hfill

\hfill

{\bf 5. Proof of Theorem 2.5, which states that the stable
bundle is hyperholomorphic if its first two Chern classes are invariant
with respect to the isotropy group of the base hyperk\"ahler
manifold. }

\hfill

 Let $B$ be a holomorphic stable bundle over  $(M,J)$. Take
the a Yang-Mills connection in $B$. Assume that $c_1(B)$ and
$c_2(B)$ are $G_M$-invariant.  Let
$\Theta$ be a curvature of $B$; we need to prove that $\Theta$ is of type
$(1,1)$ with respect to all induced complex structures on $M$, or, what is
the same, that $\Theta$ is $G_M$-invariant.

Let $I$ be an induced complex structure over $M$, such that
$I\circ J = - J\circ I$ and $I\circ J$
is another induced complex structure.
One can easily show that the Lie algebra $\g_M$ is generated
by $\mbox{\it ad}\: I$ for all such $I$. Therefore it is sufficient to
prove that $\Theta$ is of type $(1,1)$ with respect to $I$.

Let $\Theta=\Theta_{2,0}+\Theta_{1,1}+\Theta_{0,2}$ be a Hodge
decomposition of $\Theta$ with respect to $I$.
Let $\Lambda_c=\Lambda_J+\1\Lambda_K$, as in Section 4.
This is an operator of Hodge type (-2,0). For each point $x\in M$
one can define a constant

\[ Tr(\Lambda_c^2(\Theta_{2,0}\wedge\Theta_{2,0})).\]

 In the end of this section,
we will prove the analogue of Bogomolov-Gieseker ([Ko], [S])
unequality

(5.1)\hspace{6mm} $Tr(\Lambda_c^2(\Theta_{2,0}\wedge\Theta_{2,0})) < 0$

for point $x\in M$ such that $\Theta_{2,0}\neq 0$ as a section
of $\Lambda^{2,0}(End(B)$ at this point.

Let us show that

(5.2)\hspace{6mm} $\displaystyle \int_{M}
  Tr\Lambda_c^2(\Theta_{2,0}\wedge\Theta_{2,0}) =
  \int_{M} Tr \Lambda_c^2(\Theta\wedge\Theta).$

This follows from the equation

(5.3)\hspace{6mm}$\displaystyle  \Lambda_c^2(\Theta_{2,0}\wedge\Theta_{2,0}) =
  \Lambda_c^2(\Theta\wedge\Theta).$
This equatiom holds because $\Lambda_C$ maps the forms of Hodge
type $(p,q)$ into the forms of Hodge type $(p-2,q)$.
Therefore for a 4-form $\eta$ the operator $\Lambda_c^2$
acts trivially on all its Hodge components except the
(4,0)-component $\eta_{4,0}$:
\[ \Lambda^2_c(\eta)=\Lambda^2_c(\eta_{4,0}). \]
Finally, one knows that $\Theta_{2,0}\wedge\Theta_{2,0}$
is (4,0)-component of $\Theta\wedge\Theta$. This proves
equations (5.3) and (5.2).

We will show now how (5.1) and (5.2) imply
the statement of Theorem 2.5. One can deduce from the
Chern identities that for any complex bundle $B$ with a
connection and a curvature $\Theta$

(5.1A) \hspace{6mm}  $\left(2c_2(B)-\frac{r-1}{r}c_1(B)^2\right)=
    [Tr(\Theta\wedge\Theta)],$

where $[Tr(\Theta\wedge\Theta)]$ means the cohomology class
of $\C$-valued 4-form $Tr(\Theta\wedge\Theta)$ (see [Ko]).

By $r$ we denote the rank of $B$ and
by $[\alpha]$ the cohomology class in $H^*(M,\C)$
corresponding to the $\C$-valued differential form $\alpha$.

In our case $c_1(B)=0$, and therefore
\[
  \int_{M} \Lambda_c^2(\Theta\wedge\Theta) = 2\Lambda_c^2(c_2(B))
\]

This statement together with formula (5.1) imply that
$\Lambda_c^2(c_2(B))\leq 0$ and $\Lambda_c^2(c_2(B))= 0$
only if $\Theta_{2,0}=0$ everywhere. Therefore, if
$\Lambda_c^2(c_2(B))= 0$, the form $\Theta_{2,0}=0$.
Since $B$ is a Hermitian bundle,
$\Theta_{0,2}=\bar\Theta_{2,0}$. This implies that
for any induced complex structure $I$, the operator,
$ad\:\! I$ acts trivially
on $\Theta(=\Theta_{1,1})$. This implies that $\Theta$ is
$G_M$-invariant, and the bundle $B$ is hyperholomorphic.
Theorem 2.5 is proven (modulo the unequality (5.1)).

\hfill

We will proceed to the proof of the unequality (5.1).

\hfill

Since  $B$ is Yang-Mills with respect to $J$,
we have $\Lambda_J(\Theta)=0$ by Proposition 2.1.
Since the operator $\Lambda_K+\1\Lambda_I$ is of
type (2,0) with respect to $J$, we have
$\Lambda_K+\1\Lambda_I(\Theta)=0$. Since
$\Theta$ is a real form with respect to a real structure on
$\Lambda^2(End(B))$, we have also
$\Lambda_I(\Theta)=\Lambda_K(\Theta)=0$ as well.

\hfill

Let $x_1, x_1', x_2, x_2'...x_n, x_n'$ be coordinates in the bundle
$\Lambda^{1,0}_I(M)$ for some open set
$U\subset M$ such that $\bar J(x_i)=x_i'$
and $\omega_J+\1\omega_K=\sum_{i=1...n} x_i\wedge x_i'$ (where $\omega_J$ and
$\omega_K$ are K\"ahler forms associated with complex structures $J$ and $K$
respectfully). Since $\Theta$ is of type (1,1) with respect to $J$,
the operator $J$ preserves $\Theta$, i .e, $J(\Theta)=\Theta$.
There, $J(z_i\wedge z_j \wedge ...) := J(z_i)\wedge
J(z_j)\wedge ...$.
One can see that the operator $J$ maps $(2,0)$ forms in $(0,2)$ forms,
where Hodge decomposition is taken with respect to $I$.
Therefore $J(\Theta_{2,0})=\Theta_{0,2}$ and $\Theta_{2,0}$
belongs to $\Lambda\Lambda^{1,1}(End(B))$ in the sense of Section 4.

\hfill

Let us represent
\[ \Theta_{2,0}=\displaystyle\sum_{i,j=1...n}A_{ij} x_i\wedge x'_j \]
where $A_{ij}$ is a section of $End(B)$ over $U$.
This representation exists because $\Theta_{2,0}$
is a section of $\Lambda\Lambda^{1,1}(End(B))$.

Since $\bar J(\Theta_{2,0})=\Theta_{2,0}$, and
the real structure in $End (B)$ is defined by $A\arrow -A^\bot$,
we have $\sum A_{ij}x_i\wedge x_j'=\sum -A_{ij}^\bot x_i'\wedge x_j$.
Therefore,

(5.4)\hspace{10mm} $A_{ij}= A^\bot_{ij}$.

On the other hand, for any form $\sum B_{ij} x_i\wedge x_j'$,
\[ \Lambda_c\left(\sum B_{ij} x_i\wedge x_j'\right)=
   \sum_{i=1...n}B_{ii}
\] because the canonical holomorphic symplectic form is
equal to $\sum_{i=1...n} x_i\wedge x_i'$. We have proved
earlier that $\Lambda_c\left(\sum A_{ij} x_i\wedge x_j'\right)$
is a zero section of $End(B)$. Therefore

(5.5)\hspace{10mm} $\sum A_{ii}=0$
at each point of $M$.

An easy calculation shows that
\[
  \Lambda_c^2(\Theta_{2,0}\wedge\Theta_{2,0})
  =\sum_{i\neq j}-A_{ij} A{ji} + \sum_{i\neq j} A_{ii} A_{jj}.
\]

Since $Tr(AA^\bot)>0$ for every non-zero matrix $A$,
we have that

\[
   Tr\sum_{i\neq j}-A_{ij}A_{ji}
   \stackrel{(5.4)}{=\!=\!=}Tr\sum_{i\neq j}-A_{ij}A_{ij}^\bot<0
\] if any of $A_{ij},\;i\neq j$ is non-zero.

On the other hand, let us show that
$Tr \sum_{i\neq j} A_{ii}A_{jj}<0$, if not all
matrices $A_{ii},\; i=1...n$ are zero. First of all,
\[ 2\sum_{i\neq j} A_{ii}A_{jj}=(\sum A_{ii})^2 - \sum A_{ii}^2. \]

By (5.5), $\sum A_{ii}=0$ and therefore
\[ 2\sum_{i\neq j} A_{ii}A_{jj}= - \sum A_{ii}^2. \]
By the previous argument,
\[
   Tr\sum_{i=1...n}-A_{ii}^2
   \stackrel{(5.4)}{=\!=\!=}Tr\sum_{i\neq j}-A_{ii}A_{ii}^\bot\leq0
\]
Therefore

\[
  Tr(\Lambda_C^2(\Theta_{2,0}\wedge\Theta_{2,0})=
  Tr\sum_{i\neq j}-A_{ii}A_{ii^\bot}+
  Tr\sum_{i\neq j}-A_{ij}A_{ij}^\bot < 0
\]
if $\Theta_{2,0}\neq 0$. This proves unequality (5.1).$\;\blacksquare$

\hfill

\hfill

{\bf 6. Deformation spaces of holomorhpic and hyperholomorphic bundles.}

\hfill

Let $F$ be a stable bundle over a compact smooth complex
manifold $M$.

{\bf Definition 6.1} The {\bf stable deformation} $(X,x_0,{\cal F})$ of
the bundle $F$ over $M$ is the
analytic space $X$ with a marked point $x_0$ and a bundle ${\cal F}$ over
$X\times M$ such that the following holds.
The restriction of ${\cal F}$ on $\{x_0\} \times M$
is isomorphic to $F$ and the restriction of ${\cal F}$ on $\{x\} \times M$
is stable for each point $x\in X$.

For any $z\in X$ we will denote the restriction of ${\cal F}$ on
$\{z\} \times M$ by $[z]$. This is a bundle over $M$,
since  $\{x\} \times M$ is canonically isomorphic to $M$.

\hfill

{\bf Definition 6.2} Let $F$ be some sheaf over a complex manifold
$M$. Let $\cup$ be a standard $\cup$-multiplication:

\[\cup: H^1(End F)\times H^1(End F)\arrow H^2(End F\otimes End F).\]

Let $\bar k$ be a commutator map $End(F)\times End(F)\arrow End(F)$.
Obviously, $\bar k$ induces a natural map
$k:\; H^2(End F\times End F)\arrow H^2(End F)$. The composition
$k\circ\cup$ is called {\bf Yoneda pairing}.
Since $H^i(End F)= Ext^i(F,F)$, Yoneda pairing is a map

\[\iota:\;\;Ext^1(E,E)\times Ext^1(E,E)\arrow Ext^2(E,E). \]

This map is bilinear and symmetric.

\hfill

The following proposition is proven in the Section 16 of [KS].

{\bf Proposition 6.1} The Yoneda pairing on $Ext^1(F,F)$ is an
obstruction to the existing of a deformation of $F$ in the following sense.
Suppose that $\rho\in Ext^1(F,F)$ has nonzero Yoneda square. Then there
is no deformation $(D,x_0,{\cal F})$ of $F$ where $D$ is a disc in $\C$
such that the image of Kodaira-Spencer map $T_{x_0}D=\C\arrow Ext^1(F,F)$
is proportional to $\rho$. $\;\;\blacksquare$

\hfill

It turns out that in hyperholomorphic case Yoneda pairing is the only
obstruction to the existence of the deformation. We will prove this
in section 7.

\hfill

\hfill

{\bf Definition 6.3} ([G]) The marked analytic space
$Spl\galochka(F),[F]$ is called {\bf coarse moduli space of deformations}, if
there is a simple bundle $F'$ over $M$ associated to each point
$[F']\in Spl\galochka(F)$ and for any variation `
$(X, x_0, {\cal F})$ of $F$ with a
connected $X'$ there is a unique map $f:\:(X,x_0)\arrow (Spl\galochka(F),[F])$
such that for each point $x\in X$ the restriction of ${\cal F}$ on
$x\times M$ is isomorphic to the bundle, associated with $f(x)\in Spl(F)$.

Later on, we will usually omit the word coarse.

\hfill

There is a theorem of existence of a coarse moduli space of
deformations:

{\bf Theorem 6.1} (see [M], [Ko]) For any stable bundle
$F$ there exists a coarse moduli space $Spl\galochka(F)$ of the deformations
of this sheaf. If the $[F']\in Spl(F)$ is the point which corresponds
to the sheaf $[F']$ over $M$, then for the Zariski tangent space at this point
the Kodaira-Spencer map $TSpl\galochka(F)\vert_{[F']}\arrow Ext^1(F',F')$ is
isomorphism.

\hfill

Note that the space $Spl\galochka(F)$ is by no means separated or reduced.
We are chiefly interested in its reduction, denoted by $Spl(F)$.
For this space, the Kodaira-Spencer map
 $TSpl(F)\vert_{[F']}\arrow Ext^1(F',F')$ is embedding.

\hfill

The following theorem, describing
$Spl(F)$ locally, is proven in the Section 7:

{\bf Theorem 6.2} If $F$ is a hyperholomorphic
bundle, then the space $Spl(F)$ locally in the neighbourhood
of the point $[F]$
is isomorphic to an intersection of an open ball
in $Ext^1(F,F)$ with a quadratic cone
\[\{\rho\in Ext^1(F,F)\;|\;\iota(\rho,\rho)=0\}.\]

\hfill

In other words, Theorem 6.2 says that there are no obstructions
to a deformation of a hyperholomorphic bundle except the
Yoneda pairing.

\hfill

By Theorem 2.5, the stable deformation of a hyperholomorphic
bundle is also hyperholomorphic. Therefore, one can define
the hyperk\"ahler structure in the coarse deformation space
(Section 9). Unfortunately, the deformation space is
not smooth, so the definition of a hyperk\"ahler manifold
must be modified to include singular manifolds.

\hfill

First, we have to make some preliminary definitions.
Let $R$ be an $\R$-vector space with a quaternion action and a Euclidean
metric $(\cdot,\cdot)$. For each $I\in\h \;|\; I^2=-1$ the action
$I$ on $R$ defines a complex structure on $R$. Such complex
structure is called {\bf induced by quaternion action}.
The space $R$ is called {\bf quaternionic Hermitian} if
the metric on $R$ is Hermitian for any induced complex
structure $I$. For any such $I$, take $\omega_I=(\cdot,I\cdot)$.
This is a real symplectic form on $R$. Now, take
$J\in\h\;|\;I\circ J=-J\circ I,\; J^2=-1$ and $K=I\circ J$.
Of course, $J$ and $K$ define induced complex structures on $R$.
One can prove that the 2-form $\omega_J+\1\omega_K$
does not depend on the choice of $J$, $K$ as long as these quaternions
satisfy the conditions above. This form is called {\bf canonical
symplectic form} associated with $R$ and $I$.

\hfill

{\bf Definition 6.5} The (possibly singular) real analitic variety $S$
is called {\bf (singularly) hyperk\"ahler}
if it is endowed with the following structures.

\hfill

(i) There is an action of an algebra of quaternions ${\Bbb H}$ on the
real Zariski tangent space $TS$ to $S$.

\hfill

(ii) For each $\tilde I\in {\Bbb H}$ with $\tilde I^2=-1$ there is a
complex structure $I$ over $S$ such that the complex structure operator
on $TS$ associated with $I$ is equal to $\tilde I$.
This complex structure is called {\bf induced by the quaternion
action}.

\hfill

(iii) For each $x\in S$ the fiber $TS_x$ of the sheaf $TS$
in $x$ has a Euclidean metric, which defines a quaternionic
Hermitian structure on $TS_x$.

\hfill

Consider $S$ as a complex variety with an
induced complex structure $I$. Consider a
coherent sheaf $Hom_{\cal O_S}(TS\otimes TS,{\cal O}_S)$
of bilinear forms over $S$. Of course, the fiber of
$Hom_{\cal O_S}(TS\otimes TS,{\cal O}_S)$ at each
point $x\in S$ is canonically embedded in a space
of bilinear $\C$-valued forms over $TS_x$.

(iv) There is a holomorphic section $\Omega$ of
$Hom_{\cal O_S}(TS\otimes TS,{\cal O}_S)$, which,
being restricted to a fiber of $TS$ in $x\in S$,
gives a canonical symplectic
form associated with $TS_x$ and $I$.

\hfill

\hfill

Note that this definition is compatible with a usual one: a smooth
manifold is hyperk\"ahler in traditional sense if and only if it is
hyperk\"ahler in the sense of this definition.

\hfill

The following theorem is to be proven in Section 9:

{\bf Theorem 6.3.} The space $Spl(B)$ is (singularly) hyperk\"ahler for a
hyperholomorphic $B$.

\hfill

This result is a generalization of the theorem by Kobayashi:

{\bf Theorem 6.4.} ([Ko2]) Let $B$ be a simple Yang-Mills
bundle $B$ over a hyperk\"ahler manifold $M$ such that
$H^2(End(B))=\C$. Let $Spl_{ns}(B)$ be the open subspace
in $Spl(B)$ consisting of all nonsingular points $[B']$
such that $H^2(End(B'))=\C$.
The space $Spl_{ns}(B)$ has a canonical hyperk\"ahler
structure.

\hfill

In the case when $M$ is a surface, $Spl(B)$ is smooth
(Mukai [M]) and hyperk\"ahler (Mukai [M] and Itoh [I]).

\hfill


\hfill

\hfill

{\bf 7. Local deformations of a hyperholomorphic connection.
The proof Theorem 6.2, which states that the only obstruction to the
deformation of a hyperholomorphic bundle is Yoneda pairing. }

\hfill

\hfill

Let $B$ be a Hermitian bundle of degree zero
over a K\"ahler manifold $M$.

By the virtue of Uhlenbeck-Yau theorem
(see Section 1), the stable holomorphic structures
on $B$ are in the one-to-one correspondence
to the Hermitian connections $\nabla$ with
a curvature $\Theta$ such that

\hfill

(7.1)\hspace{3mm} $\Lambda(\Theta)=0$ (where $\Lambda$ is the Hodge
operator)

\hspace{-3mm}and

(7.2)\hspace{3mm} $\Theta\in \Lambda^{1,1}\otimes End(B)$, and $B$
is indecomposable with respect to $\nabla$.

\hfill

Let $B$ be the undecomposable Hermitian vector bundle with
Yang-Mills connection $\nabla$ over a K\"ahler
manifold $M$. We are interested in local deformations of
$\nabla$ which preserve conditions (7.1) and (7.2).
One can see that the space of such deformations is a
local deformation moduli space for $B$ as for a stable holomorphic bundle.

Let $\hat\rho\in\Lambda^1_R(M)\otimes{\goth u}(B)$, where ${\goth u}(B)$ is
the bundle of skew-adjoint endomorphisms of $B$ and $\Lambda^1_R$ is a
space of real differential 1-forms. Obviously, the connection $\nabla_1 =
\nabla+\hat\rho$ is Hermitian, but in general, neither Yang-Mills nor
holomorphic.

The bundle $E=End(B)={\goth u}(B)\otimes\C$ has a unique real structure
such that ${\goth u}(B)$ is a real subbundle. Therefore we can decompose
$\hat\rho$ into the sum of (1,0) and (0,1) parts:
$\hat\rho:=\rho+\bar\rho\;$ (see Section 3).

The curvature $\Theta_1$ of the connection $\nabla_1$ is equal to
$\Theta+\nabla(\hat\rho)+\hat\rho\wedge\hat\rho$.
For infinitesimal $\hat\rho$,
the $End(B)$-valued form $\hat\rho\wedge\hat\rho$ is infinitesimal of
second order compared to $\nabla(\hat\rho)$. Therefore the first step of
description of the moduli space will be determining for
which $\hat\rho$ the form
$\nabla(\hat\rho)$ satisfies conditions (7.1) and (7.2).

The following proposition partially answers this question.

\hfill

Let $B$ be a Hermitian bundle over a K\"ahler manifold $M$.

{\bf Proposition 7.1.} Let a form $\hat\rho\in\Lambda^1_R(M)\otimes{\goth
u}(B)$ be decomposed as above: $\hat\rho=\rho+\overline\rho$.
If $\rho$ is $\Delta_{\6}$-harmonic, then
$\nabla(\hat\rho)$ satisfies conditions (7.1) and (7.2).
The condition (7.1) means that $\Lambda(\nabla\hat\rho)=0$ and
(7.2) means that $\nabla\hat\rho\in\Lambda^{1,1}(M)\otimes End(B)$.

\hfill

{\bf Proof:} Proving (7.1): the form $\nabla\hat\rho$ is of type (1,1) because
its (2,0)-part $\6\rho$ and its (0,2)-part $\bar\6\bar\rho$ both vanish.

Proving (7.2): the condition (7.1) implies that
$\nabla\hat\rho=\bar\6\rho+\6\bar\rho$. Therefore, by Kodaira's
identities,
\[
   \Lambda\nabla\hat\rho=\Lambda(\bar\6\rho+\6\bar\rho)=
   \1\6^*\rho-\1\bar\6^*\bar\rho+\nabla\Lambda\hat\rho
\]
The last term vanishes because $\Lambda(\alpha)=0$ for any 1-form
$\alpha$, and $\6^*\rho=\bar\6^*\bar\rho=0$ because $\rho$
is $\6$-harmonic. Therefore $\Lambda\nabla\hat\rho=0$.

$\;\;\blacksquare$

\hfill

Note that the space of $\Delta_\6$-harmonic forms in
$\Lambda^{1,0}(M)\otimes End(B)$ is canonically conjugate to the space
$H^1(End B)\tilde=Ext^1(B,B).$

\hfill

Denote the standard Hermitian norm on a space of differential forms with
coefficients in the Hermitian bundle as $\|\cdot\|$.

Last part of this section is dedicated to the proof of the following
theorem:

{\bf Theorem 7.1} For any given hyperholomorphic bundle
$B$ with a hyperholomorphic (in particular, Yang-Mills) connection
$\nabla$ and a curvature $\Theta$ there exists a constant
$\varepsilon$ with the following property.

For each $\Delta_{\6}$-harmonic form

\[
   \rho\in\Lambda^{1,0}(M,End (B)) \;{\rm with}\;
   \|\rho\|<\varepsilon, \;\;\; \iota(\rho,\rho)=0
\]
there exist a form $\eta\in\Lambda^{1,0}(M,End(B))$ with
$\|\eta\|<1/4\;\|\rho\|$ such that a connection

\[
  \nabla_\rho =\nabla+\eta+\bar\eta+\rho+\bar\rho
\]
has a curvature
\[
  \nabla(\rho+\bar\rho)+\Theta
\]
There $\Theta$ denotes the curvature of $B$.
Thus constructed by $\rho$ form
$\eta$ holomorphically depends on  $\rho$. In the statement above
$\iota$ denotes Yoneda pairing, see Definition 2.2.

\hfill

\hfill

The Theorem 7.1 will be proven later in this section. Now we will
demonstrate some of its implications.

{\bf Proposition 7.2.} The connection $\nabla_\rho$, supplied
by Theorem 7.1, is Yang-Mills.

{\bf Proof.} Let $\hat\rho:=\rho+\bar\rho$ as in Proposition 7.1.
Since $\rho$ is $\6$-harmonic, $\6\rho=\bar\6\bar\rho=0$.
This implies that $\nabla\hat\rho$ is (1,1) form. Therefore
$\Theta+\nabla\hat\rho$ is (1,1)-form and
by Newlander-Nierenberg theorem (Theorem 2.1)
the connection $\nabla_\rho$ is holomorphic.
To see that it is Yang-Mills we need to prove that
$\Lambda(\Theta+\nabla_\rho)=0$. By our construction,
$\Lambda(\Theta)=0$, while $\Lambda(\nabla\hat\rho)=0$
by Proposition 7.1. $\;\blacksquare$

\hfill

Let us show how Theorem 7.1 imples Theorem 6.2.
Denote the intersection of a quadratic cone
\[
   \{\rho\in Ext^1(B,B)\;\;
   |\;\; \iota(\rho,\rho)=0\}
\]
with an open disc of radius $\varepsilon$ in

\[
   ker\Delta_\6\;{\bf\mid}\;_{_{\Lambda^{1,0}(M)\otimes End(B)}}
   = Ext^1(B,B)
\]
as $S$. Choose a complex structure on $S$ conjugate to the standard
one. Consider the bundle ${\cal B}$ over $S\times M$ which is, as a vector
bundle, isomorphic to a pullback of $B$ from $M$ to $S\times M$ by the
projection map. Choose a connection $\nabla_1$ on ${\cal B}$ wich is
(a) trivial in the direction of $S$ and (b)
at each point of $\rho_0\times M$ is
isomorphic to $\nabla_{\rho_0}$ in the direction of $M$.  Let $\Theta_\rho$
be a curvature of a connection $\nabla_\rho$.  Let $\Xi$ be the 2-form over
$S\times M$ with coefficients in $End B$, which is glued of
the forms $\Theta_\rho$ on the following way.
The forms $\Theta_\rho,\;\rho\in S$ are defined
over $M\times\{\rho\}$, where $M\times \{\rho\}$ is the fibre
over $\rho$ of the projection $S\times M\arrow S$.
Let $\Theta_1$ be the curvature of $\nabla_1$.

Let us decompose $\nabla$ into the sum of two components:
$\nabla=\nabla_S+\nabla_M$. The component $\nabla_S$
corresponds to the derivation in the direction of $S$,
while the component $\nabla_M$ corresponds to the derivation
in the direction of $M$. Analogously, let us decompose
the (1,0)-component $\6$ of the connection $\nabla$:
$\6=\6_S+\6_M$.

In this notation, one can easily see that
$\Theta_1-\Xi=\nabla_S\hat\rho.$

The form $\Xi$ is of type (1,1) by the construction.  Since $\eta$
holomorphically depends on $\rho$ in the standard complex structure, it
depends on $\rho$ antiholomorphically in the opposite one.  Therefore
$\6_S\eta=0$ in this structure and $\nabla_S\hat\eta$ is a form of
type (1,1). Therefore, the form $\Theta$ is of type (1,1) and the
connection $\nabla_1$ is holomorphic.

This consideration shows how Theorem 7.1 gives a construction of a
holomorphic local variation $(S,{\cal B})$ of Yang-Mills connections for a
hyperholomorphic $B$.  By the Uhlenbeck-Yau theorem, this is the same as
the construction of local variation of a stable bundle. Since
$\|\eta\|<1/4\;\|\rho\|$ and $\rho$ is represented by the classes in
$Ext^1(B,B)$, the differential of an obvious map $K:\;S\arrow Ext^1(B,B)$
is the Kodaira-Spencer map.

We have constructed a local variation of $B$ with the base $S$
and the following properties. First, the Kodaira map
$K:\;S\arrow Ext^1(B,B)$ is imbedding.
Second, its image coincides with the set of vectors
$\kappa\in Ext^1(B,B)$ such that $\iota(\kappa,\kappa)=0$.
These two properties are sufficient for the variation to
be locally universal (see [KS]), and the Theorem 6.2 is proven.

\hfill

{\bf Corollary 7.1} Theorem 7.1 yields a construction of a universal local
variation of a stable bundle $B$, if $B$ is hyperholomorphic.

\hfill

\hfill

{\bf Proof of Theorem 7.1}.

\hfill

As usually, we denote $\hat\eta:=\eta+\bar\eta$ and
$\hat\rho:=\rho+\bar\rho$.

\hfill

{\bf Lemma 7.1.} The form $\eta$ suffices the
statement of Theorem 7.1 if the following equation holds:

\hspace{5mm}(7.3)\hspace{5mm}
$-\nabla\hat\eta=\hat\rho\wedge\hat\rho+\hat\rho\wedge\hat\eta+
\hat\eta\wedge\hat\rho+\hat\eta\wedge\hat\eta$

{\bf Proof:} By definition, $\nabla_\rho=\nabla+\hat\rho+\hat\eta$.
Here, as elsewhere, the forms $\hat\rho$ and $\hat\eta$
are considered as operators as follows:
$\hat\eta(\alpha)=\hat\eta\wedge\alpha$. In this notation,
$\nabla\circ\hat\rho+\hat\rho\circ\nabla=\nabla(\rho)$.
Therefore,

\hspace{5mm}$(*)$\hspace{5mm}
$\Theta_\rho := \nabla_\rho^2 = \nabla(\hat\eta)+\nabla(\hat\rho)+
\hat\rho\wedge\hat\rho+\hat\rho\wedge\hat\eta+
\hat\eta\wedge\hat\rho+\hat\eta\wedge\hat\eta+\nabla^2$,

where $\nabla^2=\Theta$. Theorem 7.1 states that
$\Theta_\rho=\Theta+\nabla(\hat\rho)$, and by $(*)$ this
is equivalent to
$\nabla\hat\eta+\hat\rho\wedge\hat\rho+\hat\rho\wedge\hat\eta+
\hat\eta\wedge\hat\rho+\hat\eta\wedge\hat\eta=0.\;\;\blacksquare$

\hfill

We will try to find a solution of (7.3) in the form of
a Taylor serie. Let us make the following change
of notation. We shall redenote $\hat\rho$
by $t\hat\eta_1$ and $\hat\eta$ by
$(t^2\hat\eta_2+t^3\hat\eta_3+...)$.

Then the equation (7.3) could be rewritten as a system of equations:

\hfill

\hspace {10mm}$-\nabla\hat\eta_2=\hat\eta_1\wedge\hat\eta_1$

\hspace {10mm}$-\nabla\hat\eta_3=\hat\eta_1\wedge\hat\eta_2+
                                 \hat\eta_2\wedge\hat\eta_1$

(7.4) \hspace{12mm}.......

\hspace {10mm}$-\nabla\hat\eta_n=\sum_{j+j=n}\hat\eta_i\wedge\hat\eta_j$

\hspace{20mm} .......

\hfill

We will find $\eta_2$, .... $\eta_n$, .... which satisfy (7.4) and show
that for $t$ small enough the power series

\[\hat\eta=t^2\hat\eta_2+t^3\hat\eta_3+... \]

\hspace{-5mm}converges. Further on, we will as usually
denote the (1,0) component of $\hat\eta_i$
as $\eta_i$ and the (0,1) component - as $\bar\eta_i$.

Let us solve the first equation

\[
  -\nabla\hat\eta_2 = \hat\eta_1\wedge\hat\eta_1
  (= t^{2}\hat\rho\wedge\hat\rho).
\]

By our assumptions, for Yoneda pairing

\[
   \iota: H^1(End B)\times H^1(End B)\arrow H^2(End B).
\]
$\iota(\rho,\rho)=0$. One can easily prove that (see [KS])
that the cohomology class of $\rho\wedge\rho$ in $H^2(End(B))$
coincides with $\iota(\rho,\rho)$.
Therefore the form $\rho\wedge\rho$, being obviously
$\6$-closed, is also $\6$-exact.

Let $G_\6$ be the Green operator
associated with the Laplacian $\Delta_\6$, and $\Gamma=\6^*\circ G_\6$.
Clearly, $G_\6\Delta_\6\tau=\tau$ and $\Delta_\6\tau=\6\6^*\tau$
for $\6$-exact $\tau$. Since $G\6\6^*=\6\6^* G$, we see that
$\6\6^* G\tau=\tau$ for $\6$-exact $\tau$.

We have proven

{\bf Lemma 7.2.} For each $\6$-exact form
$\tau$ with coefficients in a holomorphic Hermitian
bundle over arbitrary K\"ahler manifold
the following equation holds:

\[ \6\circ\Gamma(\tau)= \tau.\;\;\blacksquare \]

\hfill

As a corollary, we obtain the following proposition:
Let $\eta_2=-\Gamma(\eta_1\wedge\eta_1)$. Then

(7.5)\hspace{10mm}$\6\eta_2=-\eta_1\wedge\eta_1$

\hspace{-5mm}and

$(7.5')$\hspace{10mm}$\bar\6\bar\eta_2=-\bar\eta_1\wedge\bar\eta_1.$

The bundle $\Lambda^1(M,End(B))$ in endowed with a
canonical real structure $T$ (see Section 3).
The operator $T$ interchanges subbundles  $\Lambda^{0,1}(M,End(B))$
and $\Lambda^{1,0}(M,End(B))$. Consider the derivation $\tilde T$
of the algebra of $End(B)$-valued differential forms which acts on
$\Lambda^1(M,End(B))$ as $T$ and on $\Lambda^i(M)\otimes End(B)$
for $i>1$ by Leibnitz formula
\[
   \tilde T(\alpha\wedge\beta)=\tilde T(\alpha)\wedge\beta+
   \alpha\wedge\tilde T(\beta).
\]

\hfill

{\bf Lemma 7.3} For any two $End(B)$-valued (1,0)-forms
$\lambda$ and $\mu$, the following holds:

\[
  \tilde T(\lambda\wedge\mu) =
  \bar\lambda\wedge\mu + \lambda\wedge\bar\mu
\]
\[ T(\6\lambda+\bar\6\bar\lambda)=2\bar\6\lambda + 2\6\bar\lambda \]

The first property is just a part of definition of $T$.
Let $\hat\lambda:=\lambda+\bar\lambda$. Let, as usualy,
$\nabla:=\6+\bar\6$. We will say that the form is real
if it is preserved by $T$. It is very easy to see that
for a real $p$-form $\eta$, $\tilde T(\eta)=p\eta$. Now,
$\tilde T(\nabla(\hat\lambda))=2 \nabla(\hat\lambda)$
because $\nabla(\hat\lambda)$ is a real form. Therefore

\[
   \tilde T(\bar\6\bar\lambda+\6\bar\lambda+
   \6\lambda+\bar\6\lambda)=
   2(\bar\6\bar\lambda+\6\bar\lambda+\6\lambda+\bar\6\lambda.)
\]

The operator $\tilde T$ interchanges $\Lambda^{1,1}(M,End(B))$
and $\Lambda^{2,0}(M,End(B))\oplus\Lambda^{0,2}(M,End(B))$.
Therefore
\[
   T(\6\lambda+\bar\6\bar\lambda)=
   2\bar\6\lambda + 2\6\bar\lambda.\;\;\blacksquare
\]

\hfill

Applying $\tilde T$ to the sum of equations (7.5)
and $(7.5')$ one obtains that

(7.6)\hspace{20mm}$\6\bar\eta_2+\bar\6\eta_2=
-\eta_1\wedge\bar\eta_1-\bar\eta_1\wedge\eta_1$

\hfill

The summation of (7.5), (7.5') and
(7.6) gives us a solution $\hat\eta_2=\eta_2+\bar\eta_2$ of an equation

\[ -\nabla\hat\eta_2 = \hat\eta_1\wedge\hat\eta_1.\]
The same consideration proves the following lemma:

{\bf Lemma 7.4} Let $\xi_1,..., \xi_n$ be (1,0)-forms with coefficients in
$End(B)$. Let $\bar\xi_i$ be conjugate forms with respect to the real
structure on $End(B)$.  Let $\hat\xi_i = \xi+\bar\xi$. Suppose that
$\sum_{i+j=n}\xi_i\wedge\xi_j$ is $\6$-exact. Let $\Gamma$ be the inverse
operator to $\6$ defined above, and
$\sigma=\Gamma(\sum_{i+j=n}\xi_i\wedge\xi_j)$.  Let
$\hat\sigma=\sigma+\bar\sigma$.  Then
$\nabla\hat\sigma=\sum_{i+j=n}\hat\xi_i\wedge\hat\xi_j$

{\bf Proof:} By our assumption,
$\displaystyle\sum_{i+j=n}\xi_i\wedge\xi_j$ is $\6$-exact.
The Lemma 7.2 implies that in this case
\[ \6\sigma=\displaystyle\sum_{i+j=n}\xi_i\wedge\xi_j.\]
Applying $T$ to this equation one obtains
\[
  \bar\6\bar\sigma=
  \displaystyle\sum_{i+j=n}\bar\xi_i\wedge\bar\xi_j,
\]
and applying $\tilde T$ to the sum of these two equations
one obtains
\[
  \bar\6\sigma+\6\bar\sigma=
  \displaystyle\sum_{i+j=n}
  (\xi_i\wedge\bar\xi_j+\bar\xi_i\wedge\xi_j).
\]

Finally, summing up these three equations one obtains

\[
   \nabla\hat\sigma=
   \displaystyle\sum_{i+j=n}\hat\xi_i\wedge\hat\xi_j.
   \;\;\blacksquare
\]

\hfill

Now, let us define the following recursive sequence:

\[\eta_1=\rho/t,\;\;\eta_2=\Gamma(\eta_1\wedge\eta_1),\;\;
...\;,\;\; \eta_n=\Gamma(\sum_{i+j=n}\eta_i\wedge\eta_j),\;\;....\]

By Lemma 7.4, if each step yields the $\6$-exact form
$\sum_{i+j=n}\eta_i\wedge\eta_j$, then

\[\nabla\hat\eta_n=\sum_{i+j=n}\hat\eta_i\wedge\hat\eta_j\]

and we obtained solutions of (7.4) we looked for.
Let us prove that the form $\tau_n
:= \sum_{i+j=n}\eta_i\wedge\eta_j$ is $\6$-exact. First of all, we show that
this form is $\6$-closed. Since

\[\6\eta_n = \sum_{i+j=n}\eta_i\wedge\eta_j\]
on each inductive step of our construction, we obtain that

\[\kappa:=\6(\sum_{i+j=n}\eta_i\wedge\eta_j)=\sum_{i+j+k=n}\eta_i\wedge\eta_j\wedge\eta_k.\]

It is easy to see that each 3-vector in $\Lambda^{3,0}$ is mapped by
$\kappa$ to the element of $End(B)$ which is equal to a sum of several
terms like $[A[B,C]]+[[B,C]A]+[[C,A]B]$. By the Jacoby identity, $\kappa$
vanishes, and the form $\sum_{i+j=n}\eta_i\wedge\eta_j$ is closed.

\hfill

To prove exactness of $\tau$ we use $\6\6^j$-lemma (Theorem 4.3).  First of
all, the form $\rho$ is $\Delta_\6$-harmonic and by Theorem 4.1 it is also
$\6^j$-closed (even $\6^j$-harmonic).
By Proposition 4.4, the operator $\6^j$ anticommutes with $\6$ and $\6^*$.
Therefore $\6^j$ commutes with $\Delta_\6$, $G_\6$ and anticommutes with
$\Gamma:=\6^*\circ G_\6$.  This implies that
each form $\eta_i$ constructed recursively above
is $\6^j$-closed. Let us prove inductively that the forms $\eta_i\;\;|\;\;
i> 1$ are $\6^j$-exact.  Suppose that $\eta_i$ is $\6^j$-exact for $1<i<n$.
The form

\[\tau_n = \sum_{i+j=n}\eta_i\wedge\eta_j \]
is $\6^j$-exact because it is a sum of several summands,
each of those is a product of two
$\6^j$-closed forms, at least one of which is $\6^j$-exact.  By Proposition
5.3, $\Gamma(\tau_n)=\eta_n$ also $\6^j$-exact.  Finally, since $\eta_n$ is
also $\6$-closed, by $\6\6^j$-lemma (Theorem 5.3)
it is $\6$-exact.

\hfill

This proves that the system (7.4) will have solutions of the following
form:

\[\eta_2=\Gamma(\eta_1\wedge\eta_1),\;\;
...\;,\;\; \eta_n=\Gamma(\sum_{i+j=n}\eta_i\wedge\eta_j),\;\;....\]

One knows that the Green operator $G$ is compact. Therefore
$\Gamma$ is also compact. Let $r$ be the norm of $\Gamma$. Obviously,

\[\|\eta_n\|\leq r\cdot n\cdot sup_{i<n}\|\eta_i\|^2\leq
n!\cdot\|\eta_1\|^{2^n}r^{2^n}= n!\cdot t^{-2^n}\|\rho\|^{2^n}r^{2^n}\]

If $\|\rho\|<\frac{\varepsilon}{4r}$ the serie

\[\eta=\sum^\infty_{i\geq 2}\eta_it^i\]

 will converge for $t=1$ and the result will be sufficiently small:
$\|\eta\|<1/4\|\rho\|$.  Finally, $\eta$ holomorphically depends on $\rho$
since $\Gamma$ is linear and the map

\[\{\eta_1....\eta_n\arrow\sum_{i+j=n}\eta_i\wedge\eta_j\}\]

is holomorphic. Theorem 7.1 is proven.$\; \; \blacksquare$

\hfill

\hfill

{\bf Section 8. Comparing Laplacians.}


\hfill


Let $B$ be a hyperholomorphic bundle over a hyperk\"ahler manifold
$M$. Let $\Delta_{\6_{L}}$ be the Laplace operator $\Delta_\6$
associated with arbitrary induced complex structure $L$.

{\bf Theorem 8.1} If $I$, $L$ are two induced complex structures
over $M$, we have an identity:

\[ (\Delta_{\6_{I}})^L=(\Delta_{\6_{I^{{}^L}}}). \]

\hfill

We use there the notation of the Section 4, where $(\:)^L$
denotes $L\circ(\:)\circ L^{-1}$.

{\bf Proof:} Let $d$ denote the connection operator in $B$
(denoted by $\nabla$ elsewhere). Proposition 3.1 immediately
implies that

\[ \Delta_{\6_{I}}= 1/2(\Delta_d + \1 [\Lambda_I,d^2]). \]

Since $B$ is hyperholomorphic, its curvature operator
$d^2$ is $L$-invariant with respect to an induced complex structure
operator $L$. Therefore $(d^2)^L=d^2$. By Proposition 3.1(c),
$\Delta_d^L=\Delta_d$. Therefore

\[ (\Delta_{\6_{I}})^L=1/2 (\Delta_d^L +
\1[\Lambda^L_I,(d^2)^L])=1/2(\Delta_d+\1[\Lambda^L_I,d^2]. \]

Finally, the identity $\Lambda_I^L=\Lambda_{L\circ I\circ L^{-1}}$
is proven by a simple computation with quaternions. $\blacksquare$

\hfill

Theorem 8.1 implies that the group of unitary quaternions acts
transitively on the set of Laplace operators $\Delta_{\6_{L}}$
where $L$ runs through the set of induced complex structures
(see Lemma 8.1). Let $H^i_L(B)$ denote the space of
holomorphic cohomology of a hyperholomorphic bundle $B$,
taken with respect to the induced complex structure $L$.

{\bf Corollary 8.1:} The space $H^i_I(B)$ does not depend on the
choice of an induced complex structure $I$.

{\bf Proof:} The space $H^i_I(B)$ (to be precise, its
complex conjugate) can be canonically identified with the space of
$\Delta_{\6_{I}}$-harmonic sections of $\Lambda_I^{i,0}(B)$,
where $\Lambda_I^{i,0}$ means that Hodge grading is taken
with respect to $I$. One can easily see that the action of
$L$ on $\Lambda^1(M,B)$ produces an isomorphism between
$\Lambda_I^{i,0}(B)\subset\Lambda^1(M,B)$
and  $\Lambda_{I^{{}^L}}^{i,0}(B)\subset\Lambda^1(M,B)$. By Theorem 8.1,
under this isomorphism $\Delta_{\6_{I}}$ goes into
$\Delta_{\6_{I^L}}$. Therefore $L$ maps $H^i_I(B)$
into $H^i_{I^L}(B)$. Finally, Corollary 8.1
is implied by the following lemma:

{\bf Lemma 8.1:} For each pair of induced complex structures
$L$ and $L'$ over $M$ there is an induced $R(L,L')$ such that
$L^{R(L,L')}=L'$.

{\bf Proof:} This is an easy linear algebraic
computation.

Let $G_M$ be the isotropy group for $M$.
This group is isomorphic to $SU(2)$.
There is a unique element $r\neq 1$
in $SU(2)$ such that $r^2=1$. This element is defined by
the matrix $-1$. The set of induced complex structures is isomorphic
to the set of elements $a\in G_M$ such that $a^2=r$.

The group $G_M$ is isomorphic to the group of
unitary quaternions $\{h\in\h\:|\:h\bar h=1\}$.
The element $r$ is just $-1\in \h$. The unitary quaternion
$s$ has the square $-1$ if $s=aI+bJ+cK$ for
real $a,b,c$ such that $a^2+b^2+c^2=1$.
Therefore the set $S$ of such quaternions is a unit sphere in
$\R^3$. Clearly,

\[ I\circ (aI+bJ+cK)\circ I^{-1}= aI-bJ-cK \]

This implies that for induced complex structure
$B\in S$, the map $A\arrow A^B$ from $S$ to $S$ is
a symmetry ($180^\circ$ angle rotation)
of $S$ around the line going through $B$
and $-B$. Now, for each $L$, $L'$ in $S$ there is
an element $R(L,L')$ in $S$ such that the symmetry
around the line which goes through $R(L,L')$
and $-R(L,L')$ maps $L$ into $L'$. Simply,
connect $L$ and $L'$ with a big circle. Take the point
on this circle which equally divides the interval between $L$ and
$L'$. This is $R(L,L')$.

 The Lemma 8.1, and therefore Corollary 8.1 are proven.
$\blacksquare$

\hfill

By $Ext^i_I(B,B)$ we mean $Ext^i(B,B)$ where $B$ is taken
as a holomorphic bundle over $(M,I)$.

{\bf Corollary 8.2:} The graded rings $\oplus_i Ext^i_I(B,B)$ are isomorphic
for all induced $I$ over $M$.

{\bf Proof:} For each $i$, $Ext^i_I(B,B)\cong H^*_I(End(B)$, and
$End(B)$ is hyperholomorphic for hyperholomorphic $B$. The isomorphism
constructed in Corollary 8.1 is clearly multiplicative.
$\:\;\; \blacksquare$

\hfill

\hfill

{\bf 9. The space of deformations of a hyperholomorphic bundle
is hyperk\"ahler (Theorem 6.3).}

\hfill

This proof is completely analogous to the proof
of the same theorem for the case when $M$ is a surface ([M])
and for the case when $H^2(End(B))= \C$ (Theorem 6.4,
proven in [Ko2]).

We need to show that $Spl(B)$ suffices
points (i)-(iv) of Definition 6.5.

\hfill

(i)-(ii):
Corollary 8.2 together with Theorem 6.2 provide that the local
deformation space $(S,I)$ of a hyperholomorphic bundle $B$ over
$(M,I)$ is independent on the choice of an induced complex
structure $I$. Since the $(S,I)=:S$ is a subset of an open ball
in $Ext^1(B,B)$, one can easily check the following, using Corollary
8.2. Let $\tilde I$, $\tilde J$ and $\tilde K$ are complex structures
on $S$ defined by induced complex structures $I$, $J$, $K$ over $M$.
Assume that $I\circ J= - J\circ I=K$. Then
$\tilde I\circ\tilde J= -\tilde J\circ\tilde I=\tilde K$.

The last statement is sufficient to prove that $Spl(B)$
satisfies (i)-(ii) from the Definition 6.5.

The quaternion action on the tangent space $T_{[B']}$
to $Spl(B)$ in the point, corresponding to the bundle $B'$,
is compatible with the quaternion action on $H^1(End(B'))$
constructed in Corollary 4.3.

\hfill

(iii) Since
$T_{[B]} Spl(B)$ is a subspace of $H^1_I(End(B))$, one
can restrict the Hermitian metric on
$H^1_I(End(B))$ to $T_{[B]} Spl(B)$. Now, the Hermitian metric
on $H^1_I(End(B))$ is defined as follows. Cohomology
classes from $H^1_I(End(B))$ can be canonically identified
with $\Delta_\6$-harmonical forms in $\Lambda^{1,0}(End(B))$.
The bundle $\Lambda^{1,0}(End(B))$ is Hermitian,
so there is an Hermitian metric on the space
of its continous sections. To find the Hermitian product of two continous
sections of a Hermitian bundle over a compact
oriented Riemann manifold $M$, one takes their product
pointwise, and them integrates over $M$ the resulting
$\C$-valued function times volume form of $M$.

Take the real part of this Hermitian metric.
For induced $L$ and $End(B)$-differential forms
$\eta_1$ and $\eta_2$ we have
$(\eta_1,\eta_2)=(\eta_1^L,\eta_2^L)$,
where $(\cdot,\cdot)$ means the point-wise scalar
product. Combining this with the proof of Corollary 8.1, we
obtain that this positively defined
$\R$-valued scalar product is independent on the choice of $I$.
Since it is Hermitian with respect to induced complex
structures by the construction, this metric is quaternionic
Hermitian. The part (iii) of Definition 6.5 is proven.

\hfill

(iv)
Fix an induced complex structure $I$ on $M$.
 Let $S$ be the neighbourhood of $[B]\in Spl(B)$,
constructed in Corollary 7.1. The space $M\times S$,
according to the construction preceding Corollary 7.1,
is endowed with the classifying bundle ${\cal B}$
with the following property. Take a hyperholomorphic
bundle $B'$ over $M$ such that $[B']\in S\subset Spl(B)$.
The restriction of ${\cal B}$ on $M\times [B']\subset M\times S$
is isomorphic to $B'$.

Consider the projection $\pi: M\times S\arrow S$ and take
$E_1:=R^1\pi_*(End({\cal B}))$. One can immediately see that
the fiber of $E_1$ at each
point $[B']\in S$ is isomorphic to $H^1(End(B))=Ext^1(B,B)$.
Consider the cup-product

\[
    \cup:\;\; R^1\pi_*(End({\cal B}))\times
    R^1\pi_*(End({\cal B}))\arrow R^2\pi_*(End({\cal B}))
\]
defined by the bilinear product on $End({\cal B})$.
This is a bilinear map of coherent sheaves over $S$.
Denote $R^2\pi_*(End({\cal B}))$ by $E_2$. The cup-product map,
in this notation, is a bilinear map $\cup$ from $E_1\times E_1$
to $E_2$, or a coherent sheaves' map from
$E_1\otimes_{{\cal O}_S}E_1$ to $E_1$. By Theorem 6.2,
the Zariski tangent sheaf $TS$ is imbedded in $E_1$.
Moreover, $TS$ as a coherent subsheaf of $E_1$ is generated by all
sections $\rho$ of $E_1$, such that $\cup(\rho,\rho)=0$.

There is a trace map from $End({\cal B})$ to ${\cal O}_{M\times S}$.
By functoriality, this map defines the map

\[ \rho:\;\; R^2\pi_*(End({\cal B}))\arrow R^2\pi_*({\cal O}_{M\times S}).\]

{}From projection formula, we see that
that $R^2\pi_*({\cal O}_{M\times S})$
is a trivial sheaf over $S$ with a fiber
$H^2(M,{\cal O}_M)$. The space $H^2(M,{\cal O}_M)$
for a hyperk\"ahler $M$ is endowed with a canonical
linear form, which is defined by a canonical symplectic
form over $M$. This form is equal to the operator
$\Lambda_c$ from $H^2(M,{\cal O}_M)$ to
$H^0(M,{\cal O}_M)=\C$ (see Theorem 4.2).
Therefore there is a canonical map

\[
   \tau: \;\;R^2\pi_*({\cal O}_{M\times S})
   \arrow {\cal O}_S
\]

Composing $\cup: \;E_1\otimes_{{\cal O}_S} E_1\arrow E_2$
with $\rho:\;E_2\arrow R^2\pi_*({\cal O}_{M\times S})$ and
$\tau: R^2\pi_*({\cal O}_{M\times S})\arrow {\cal O}_S $,
one obtains the form $\tilde\Omega$ on $E_1$. Since
$TS\subset E_1$, one can restrict this form to $TS$ to
obtain the holomorphic section $\Omega_S$ of
$Hom(TS\otimes_{{\cal O}_S} TS, {\cal O}_S)$.

\hfill

We shall prove that this section is one required by (iv)
of Definition 6.5.

\hfill

Consider the restriction of $\Omega_S$ on the fiber $TS_{x}$
of $TS$ in $x\in S$, where $x$ corresponds to a hyperholomorphic bundle
$B$ over $M$. By Theorem 6.2, $TS_{x}$ is a subspace of
$Ext^1(B,B)=H^1(End(B)$ spanned by all $\gamma$ with
$\iota(\gamma,\gamma)=0$. There is a trace map
$\tilde\rho:\: H^2(End(B)\arrow H^2({\cal O}_M)$,
the map $\tau=\Lambda_c$ from $H^2({\cal O}_M)$ to $\C$
and the cup-product $\cup:\;H^1(End(B))\times H^1(End(B))\arrow
H^2(End(B))$.

The restriction $\Omega$ of $\Omega_S$ to $TS_x\otimes TS_x$
is a composition of embedding
$i: TS_x\otimes TS_x \hookrightarrow H^1(End(B))\otimes H^1(End(B))$,
then the cup-product, then the trace $\tilde\rho$ and then $\tilde\tau$:
\[
   \Omega:= \Omega_S|{{}_{_{TS_x\otimes TS_x}}}=
   i\circ\cup\circ\tilde\rho\circ\tilde\tau:
   TS_x\otimes TS_x \arrow \C.
\]

{\bf Proposition 9.1:} The form $\Omega$ is equal to the canonical
symplectic form constructed by the complex structure $I$ on the
quaternionic Hermitian space $TS_x$.

{\bf Proof:} First of all, one can redescribe $\Omega$ in terms of
harmonic forms as follows. Take two classes $[\alpha]$ and $[\beta]$
in $TS_x\subset H^1(End(B)$ and let $\alpha$, $\beta$ be the unique
$\6$-harmonic (1,0)-forms representing these classes. Obviously,

\[
   \Omega([\alpha],[\beta]) =
   \int_M Tr \Lambda_c(\alpha\wedge\beta)dVol(M).
\]

\hfill

There is an action of quaternions in $\Lambda^{1,0}(End(B)$,
which commutes with Laplacian and
induces the constructed above quaternion
action on $H^1(End(B))$ (see Corollary 4.3).
This action is Hermitian with respect to the metric on
$\Lambda^{1,0}(End(B)$ arising from a metric on $B$ and
on $\Lambda^{1,0}(M)$ which is identified with $\Lambda^1(M,\R)$
as a real bundle. The canonical symplectic form on the
fibres of this bundle is given by the formula

\[
   \gamma_1,\gamma_2\arrow Tr \Lambda_c(\gamma_1\wedge\gamma_2).
\]
This statement immediately follows from the definition of
Hermitian form in $\Lambda^{1,0}(End(B)$, definition of $\Lambda_c$
and the definition of a canonical symplectic form.
Since the Hermitian product on the
space of $\6$-harmonic forms is given by the integration of
product in fibres along $M$, one sees that the canonical symplectic
form on $TS_x\subset H^1(End(B)$ is given by the formula
\[
   \gamma_1,\gamma_2\arrow
   \int_M Tr \Lambda_c(\gamma_1\wedge\gamma_2) dVol(M).
\]
This is exactly the same formula as one for the
restriction of $\Omega_S$ to $TS_x$. Proposition 9.1
and Theorem 6.3 are proven. $\;\;\blacksquare$

\hfill

\hfill

{\bf 10. Some applications.}

Let $I_1$ and $I_2$ be two complex structures induced
by the same hyperk\"ahler structure. Let $B$ be a
bundle with connection. Assume that this connection is
hyperholomorphic with respect to the hyperk\"ahler structure.
This connection is integrable
over $(M,I_1)$ and $(M,I_2)$. Consider the holomorphic
bundles $B_1$ over $(M,I_1)$ and $B_2$ over
$(M,I_2)$ which are constructed by integrating
this connection in $B$. Theorem 6.3 endows the
manifolds $Spl(B_1; (M,I_1))$
and $Spl(B_2; (M,I_2))$ with a
hyperk\"ahler structure.

{\bf Corollary 10.1} Two manifolds
$Spl(B_1; (M,I_1))$ and $Spl(B_2; (M,I_2))$
are isomorphic as hyperk\"ahler manifolds. The complex structures on \\
$Spl(B):= Spl(B_i; (M,I_i))$
which come from $(M,I_1)$ or from $(M,I_2)$ are
induced by this hyperk\"ahler structure.

$\blacksquare$

\hfill

{}From now till the end of the section we will suppose that the
first Chern class of $B$ vanishes: $c_1(B)=0$.

\hfill

{\bf Definition 10.1} The holomorphic bundle $B$ of
degree zero is called {\bf strongly simple} if it has no subsheaves of
degree zero.

\hfill

{\bf Definition 10.2} The complex structure on
a K\"ahler manifold $M$ is
called {\bf of general type}
if $H^{11}(M)\cap H^2(M,\Z)=\{0\}$.

\hfill

{\bf Definition 10.3:} The hyperk\"ahler structure ${\cal H}$ is called
{\bf of general type} if there is a complex structure $I$
of general type on $M$ such that $I$ is induced by ${\cal H}$.

\hfill

{\bf Proposition 10.1} If a complex structure on $M$ is of general type
then the strong simplicity of a bundle $B$ over $M$ is equivalent
to its stability.

{\bf Proof} (see also [V], Proposition 3.1)
We know that for any coherent sheaf $F$ over a K\"ahler manifold
$M$, the cohomology class of $c_1(F)$ lies in
group $H^{11}(M)\cup H^2(M,\Z)$. Since $M$ is of general type,
this group is zero. Therefore $c_1(F)=0$ and $deg(F)=0$
for any sheaf $F$ over $M$.

By definition of stability, this means that the bundle $B$ is
stable iff it has no subsheaves $G\subset B$ with $rank(G)<rank(B)$.
$\blacksquare$

\hfill

\hfill

The following proposition is an easy consequence of Calabi-Yau
theorem ([Y]) and it was proven in [Tod].

{\bf Proposition 10.2} The set ${\cal H}$ of hyperk\"ahler structures which
induce the given complex structure $I$ on $M$ is canonically isomorphic
to a convex cone of signature $(+,-,-,...,-)$
in $H^{1,1}(M)\cap H^2(M,\R)$. The set of hyperk\"ahler structures of
general type is Zariski dense in ${\cal H}$.

\hfill

We call two bundles $B_1$ and $B_2$ algebraically equivalent if
there is a connected variety $S$ and a bundle ${\cal B}$ over $M\times S$
such that ${\cal B}$ restricted to $M\times \{s_1\}$ is isomorphic
to $B_1$ and  ${\cal B}$ restricted to $M\times \{s_2\}$ is isomorphic
to $B_2$, where $s_1$ and $s_2$ are points in $S$.

Let $K_M$ be the set of equivalence classes of hyperholomorphic
bundles with vanishing $c_1$ modulo algebraic equivalence.
Let $\{B_i,\:i\in K_M\}$ be the set of representatives of bundles in $K$,
one in each component. Denote the disjoint union of all spaces
$\{Spl(B_i),\; i\in K_M\}$ by ${\cal M}_M$.

The following proposition was proven in [V] for K3 or
abelian surfaces:

{\bf Proposition 10.3} (See also Corollary 3.1 from [V])
Suppose that $M_1$ and $M_2$ are hyperk\"ahler
manifolds of general type which belong to the same deformation
class. The spaces ${\cal M}_{M_1}$ and  ${\cal M}_{M_2}$
are diffeomorphic.

\hfill

In other words, the classifying space of hyperholomorphic
bundles over the generic hyperk\"ahler $M$
does not depend on the deformation class of $M$.

{\bf Proof:} We call the following operation on a complex
manifold $(M,I)$ {\bf standard}.

a) choose some hyperk\"ahler structure ${\cal H}$ of general type
 such that $I$ is induced by ${\cal H}$.

b) choose another complex structure $L$ which is induced by ${\cal H}$.

The result of this operation will be the complex manifold $(M,L)$.

\hfill

By the Corollary 10.1 and Proposition 10.1, the standard operation does not
change the diffeomorphism class of ${\cal M}_{(M,I)}$
if $(M,I)$ was of general type. Therefore
to prove Proposition 10.3 we should prove the following
theorem. For a surface, this is Theorem 4.1 from [V],
for arbitrary dimension see [Tod]:

{\bf Theorem 10.1} Let $M_1$ and $M_2$ be two hyperk\"ahler
manifolds with a hyperk\"ahler structure of general
type. Suppose that $M_1$ and $M_2$ belong to the same deformation
class. Then $M_1$ can be obtained from $M_2$ by the sequence of
standard operations. $\blacksquare$

\hfill

\hfill

{\bf 11. Projectively hyperholomorphic bundles. }

\hfill

Most facts about hyperholomorphic budles with connection
can be generalized to the case of projectively hyperholomorphic
bundles. As usually, $M$ is a compact hyperk\"ahler manifold,
and $(M,I)$ is the same manifold, considered as a complex one
with induced complex structure $I$.

\hfill

{\bf Definition 11.1.} Let $B$ be a bundle of rank $r$ with
connection and a curvature $\Theta\in \Lambda^2(End(B))$.
Take  the 2-form $Tr(\Theta)\in \Lambda^2(\C)$. The
$End(B)$-valued 2-form

\[ \Theta_{tl}:= \Theta-\frac{1}{r} Tr(\Theta)\]
is called traceless curvature of $B$.

{\bf Definition 11.2.}  Let $B$ be a Hermitian holomorphic bundle over $(M,I)$
with a traceless curvature $\Theta_{tl}$. The bundle $B$ is called
projectively hyperholomorphic if $\Theta_{tl}$ is $G_M$-invariant.
By $G_M$ we mean the isotropy group of $M$, see Definition 1.2.

\hfill

{\bf Proposition 11.1.} The Hermitian holomorphic bundle $B$ over $M$ is
projectively hyperholomorphic if and only if $End(B)$ is
hyperholomorphic as a bundle with connection induced from $B$.

{\bf Proof} (See [Ko] for analogous results:)
The curvature $\Theta_1$ of $End(B)$
maps a $End(B)$-valued form $\alpha$ to $\alpha\wedge\Theta
-\Theta\wedge\alpha$, where $\Theta$ is a curvature form in $B$.
This implies $\Theta_1\alpha= \theta\wedge\alpha -\alpha\wedge\Theta=
\Theta_{tl}\wedge\alpha - \alpha\wedge\Theta_{tl}$, because
the scalar 2-form $Tr(\Theta)$ commutes with $\alpha$.
The last form is  $G_M$-invariant if $B$ is projectively
hyperholomorphic, and therefore $End(B)$ is hyperholomorphic
if $B$ is projectively hyperholomorphic.

Conversely, assume that $\Theta_1$ is $G_M$-invariant, or,
what is the same, that $End(B)$ is hyperholomorphic.
Let $\Theta=\sum_i A_i\omega_i$, where $\omega_i\in \Lambda^2(M)$
and $A_i\in End(B)$. Take $A\in\Gamma(End(B))$ as
the section of $\Lambda^0(M,End(B))$. Since $A$ and $\Theta_1$ are
$G_M$-invariant differential forms, the form $\Theta_1(A)=
\sum_i[A,A_i]\omega_i$ is also $G_M$-invariant. Therefore
for any non-$G_M$-invariant $\omega_i$ the commutator
$[A,A_i]=0$ for any $A$. Therefore $A_i=f\cdot Id$ for
such $i$, where $f$ is a scalar function, and
$\Theta_{tl}$ is $G_M$-invariant. $\;\;\blacksquare$

\hfill

\hfill

Now, let $B$ be Yang-Mills bundle over $M$. One can apply
our version of Bogomolov-Gieseker unequality ((5.1) and (5.1A))
to $End(B)$ considered as a (not necessarily holomorphic)
bundle with connection over $(M,J)$ and obtain the following result.

{\bf Theorem 11.1} The Yang-Mills bundle $B$ of rank $r$
over $(M,I)$ is projectively hyperholomorphic if and
only if the cohomology class

\[ c_2(B)-\frac{r-1}{2r}c_1(B)^2 \]
is $G_M$-invariant.

\hfill

This also follows from Theorem 2.5 applied to the bundle
$End(B)$ and Proposition 11.1.

\hfill

One can apply Proposition 11.1 to prove the following:

{\bf Proposition 11.2:} If $B$ is Yang-Mills stable bundle
over a hyperk\"ahler surface (K3 or abelian surface),
then $B$ is projectively hyperholomorphic.

{\bf Proof:} This follows from Theorem 2.4,
and Proposition 11.1. For stable Yang-Mills $B$ the bundle
$End(B)$ satisfy the following. If $\Theta$ is its curvature,
then $\Lambda(\Theta)=0$. Applying the proof of Theorem
2.4 to $End(B)$, one sees that $End(B)$ is hyperholomorphic.
$\:\;\blacksquare$

\hfill

Since Theorem 6.2 depends only on $\6\bar\6$-lemma
in cohomology of $End(B)$ (Theorem 4.3), one immediately
obtains its generalization to projectively hyperholomorphic case:

{\bf Theorem 11.2:} If $B$ is a projectively hyperholomorphic
Yang-Mills bundle, the space $Spl(B)$ locally around the
point $[B]$ is isomorphic to the intersection of an open ball
in $Ext^1(B,B)$ with the quadratic cone
$\{\rho\in Ext^1(F,F)\;|\;\iota(\rho,\rho)=0\}$.

\hfill

\hfill

To prove the generalization of Theorem 6.3 to the projectively
hyperholomorphic bundles, one does the following.
Assume that $B$ is projectively hyperholomorphic. Proposition 4.3
supplies the construction of quaternionic action on
$T_{[B]}Spl(B)\hookrightarrow H^1(End(B))$. This shows that
$S=Spl(B)$ suffices (i) of Definition 6.5. The quaternionic
Hermitian metric on $H^1(End(B))$
((iii) of Definition 6.5) is being constructed
exactly as in hyperholomorphic case. Moreover, the holomorphic
section of $Hom(TS\otimes TS,{\cal O}_S)$ is constructed
as in the proof of Theorem 6.3 (Section 9), and one can
immedeately see that $Spl(B)$ suffices (iv) of Definition
6.5.

The only possible difficulty one meets is to prove that
the for any $L\in \h \;|\; L^2=-1$ the action of $L$
on $TS$ is integrable, i. e., is induced by some complex
structure on the variety $S$, considered as the real analytic space.
To prove that, one should realize the upen neighbourhood
$U$ of $[B]$ as an intersection of of the open
ball $H^1(End(B))$ and the quadratic cone $\iota(\eta,\eta)=0$
as in Theorem 11.2. This cone is quaternionic invariant,
so one can induce the quaternionic action on $U$ from
$H^1(End(B))$. This latter action is obviously integrable.
Moreover, as one can prove (see [Ko2]), this action
is compatible with constructed above canonical
symplectic form. This implies that this particular
quaternionic action (a priori, dependent on $B$)
coincides with the canonical quaternionic action
constructed above.

We have proven the following theorem:

{\bf Theorem 11.3:} For a projectively hyperholomorphic
stable bundle $B$ over a hyperk\"ahlem manifold $M$, the space
$Spl(B)$ is (singularly) hyperk\"ahler.

\hfill

\hfill

{\sf Acknowledgements:} I am very grateful to my advisor David
Kazhdan for a warm support and encouragement. Many thanks
also due to Michael Finkelberg for interesting discussions.
I want also to express my gratitude to Alexandr Beilinson, who explained
me Simpson's results, Victor Ginzburg for giving pointers
to Deligne, Griffiths, Morgan and Sullivan paper about
vanishing of Massey product on K\"ahler manifolds
and $\6\bar\6$ lemma, and Arkady Vaintrob
who introduced me to hyperk\"ahler manifolds.

\hfill

\hfill

\centerline{\bf Reference:}

\hfill

\hfill

[AK] Altman A., Kleiman S. Compactifying the Picard scheme //
Adv. in Math 35, p. 50-112 (1980)

\hfill

[B] Beauville A. Varietes K\"ahleriennes dont la pere classe de Chern est
nulle. // J. Diff. Geom. 18, p. 755-782 (1983).

\hfill

[Bes] Besse, A. Einstein Manifolds // Springer-Verlag, New York (1987)



\hfill


\hfill

[GH] Griffiths Ph. and Harris J. Principles of algebraic geometry. //
Wiley-Interscience, New York (1978).

\hfill

[G] Gieseker D. The moduli of vector bundles on an algebraic surface //
Ann. of Math., 106 p. 45-60 (1971).

\hfill

[I] Itoh, Mitsuhiro, Quaternion structure on the moduli space
of Yang-Mills connections  // Math. Ann. 276(1986/1987) 581-593

\hfill


\hfill

[K] Kobayashi S. Differential geometry of complex vector bundles. //
Princeton University Press, 1987.

\hfill

[Ko2] Kobayashi S. Simple vector bundles over symplectic K\"ahler
manifolds. // Proc. Japan Acad. 62(1986) 21-24

\hfill

[KS] Kodaira K., Spencer D. C. On deformations of complex structures II
// Kodaira K., Collected Works vol. II, Princeton Univ Press (1975).

\hfill

[M] Mukai S. Symplectic structure of the moduli space of sheaves on
abelian or K3 surfaces. // Invent. Math. 77 p. 101-116 (1984).

\hfill

[S] Simpson, C. T. Journal of Amer. Math. Soc. //
vol 1(4) p. 867-918. (1988)

\hfill

[T] Todorov A., Moduli of Hyper-K\"ahlerian manifolds I,II. // Preprint
MPI (1990)

\hfill

[Ty] Tyurin A., Algebraic geometric aspects of smooth structure I:
the Donaldson polynomials. // Russian Math. Surveys 44(1989).

\hfill

[Ty1] Tyurin A., The Weil-Peterson metric on the moduli space of
stable bundles and sheaves on an algebraic surface. //
Math USSR Izvestiya Vol. 38(3) p 599-621 (1992).

\hfill

[UY] Uhlenbeck K. and Yau S. T. On the existence of Hermitian Yang-Mills
connections in stable vector bundles // Comm. on Pure and Appl. Math.,
39, p. S257-S293 (1986).

\hfill

[V] Verbitsky M. On the stable bundles with vanishing $c_1$ over K3 and
abelian surface and hyperk\"ahler structures. // preprint, 1991.

\hfill

[V1] Verbitsky M. On the action of a Lie algebra SO(5) on the cohomology
of a hyperk\"ahler manifold. // Func. Analysis and Appl. 24(2)
p 70-71 (1990).

\hfill

[Y] Yau, S. T. On the Ricci curvature of a compact K\"ahler manifold
and the complex Monge-Amp\`ere equation I. // Comm. on Pure and Appl.
Math. 31, 339-411 (1978).

\end{document}